\documentclass[11pt]{article}
  \usepackage{amsthm}
  \usepackage{amssymb}
  \usepackage{amsmath}
  \usepackage{mathrsfs}
  \usepackage[all]{xy}
  \usepackage{bbm}

\newcommand{\g}{\mathrm{g}} 
\newcommand{\C}{\mathcal{C}} 
\newcommand{\dc}{\mathcal{O}} 
\newcommand{\si}{\sigma}
\newcommand{\Si}{\Sigma}

\newcommand{\ga}{\gamma}
\newcommand{\be}{\beta}

\newcommand{\io}{\iota}
\newcommand{\al}{\alpha}

\newcommand{\f}{\mathrm{f}}

\newcommand{\sst}[1]{\scriptscriptstyle{#1}}

\newcommand {\defi}{\equiv}

\newcommand{\Uc}{\mathcal{U}} 
\newcommand{\fr}{\mathrm{f}} 
\newcommand{\dr}{\mathrm{d}} 
\newcommand{\gr}{\mathrm{g}} 
\newcommand{\sr}{\mathrm{s}} 
\newcommand{\ur}{\mathrm{u}} %
\newcommand{\Urm}{\mathrm{U}} 
\newcommand{\Zrm}{\mathrm{Z}} 
\newcommand{\Crm}{\mathrm{C}} 
\newcommand{\Brm}{\mathrm{B}} 
 \author{John E. Roberts and  Giuseppe Ruzzi \\ \\
  Dipartimento di Matematica, \\
  Universit\`a di Roma ``Tor Vergata'' \\
     Via della Ricerca Scientifica  1 \\
     00133, Roma,  Italy  \\ \\
  \small{ \texttt{roberts@mat.uniroma2.it,  ruzzi@mat.uniroma2.it}}}

\title{A cohomological description of connections and curvature 
       over posets}
\date{}
 
\begin{document}

  \maketitle

\begin{abstract}
What remains of a geometrical notion 
like that of a  principal bundle when  the base space 
is not a manifold but a coarse graining of it, like 
the poset formed by a  base for the 
topology ordered under inclusion? 
Motivated  by finding a geometrical framework 
for developing gauge theories  in algebraic quantum field theory, 
we give, in the present paper, a first  answer to this question. 
The notions of transition function, connection form 
and  curvature form   find  a nice description 
in terms of cohomology, in general non-Abelian, of a poset 
with values in a group $G$.
Interpreting a $1$--cocycle as a principal bundle,  
a  connection turns out to be a  $1$--cochain associated in a 
suitable way with  this $1$--cocycle; the curvature  
of a connection turns out to be its $2$--coboundary.
We show the existence of nonflat connections, and relate 
flat connections  to homomorphisms of the fundamental group of the poset 
into $G$.  We discuss holonomy and prove  an analogue  of the
Ambrose-Singer theorem. 
\end{abstract}

  \theoremstyle{plain}
  \newtheorem{df}{Definition}[section]
  \newtheorem{teo}[df]{Theorem}
  \newtheorem{prop}[df]{Proposition}
  \newtheorem{cor}[df]{Corollary}
  \newtheorem{lemma}[df]{Lemma}

  \theoremstyle{definition}
  \newtheorem{oss}[df]{Remark}
  \newtheorem{ex}[df]{Example}

\theoremstyle{definition}
  \newtheorem{ass}{\underline{\textit{Assumption}}}[section]



\section{Introduction}

  One of the outstanding problems of quantum field theory is to
characterize gauge theories in terms of their structural properties.
Naturally, as gauge theories have been successful in describing
elementary particle physics, there is a notion of a gauge theory
in the framework of renormalized perturbation theory. Again,
looking at theories on the lattice, there is a well defined
notion of a lattice gauge theory.\smallskip

  This paper is a first step towards a formalism which adapts the
basic notions of gauge theories to the exigencies of algebraic quantum
field theory. If successful, this should allow one to uncover structural
features of gauge theories. Some earlier ideas in this direction
may be found in \cite{Rob3}.\smallskip

  In mathematics, a gauge theory may be understood as a principal
bundle over a manifold together with its associated vector bundles.
For applications to physics, the manifold in question is spacetime
but, in quantum field theory, spacetime does not enter directly as
a differential manifold or even as a topological space. Instead,
a suitable base for the topology of spacetime is considered as
a partially ordered set ({\it poset}), ordered under inclusion. This
feature has to be taken into account to have a variant of gauge theories
within algebraic quantum field theory. To do this we adopt a cohomological
approach. After all, a principal fibre bundle can be described in terms of
its transition functions and 
these form a $1$--cocycle in \v Cech cohomology with values in a group $G$. 
We develop here a $1$--cohomology
of a poset with values in $G$ and regard this as describing principal 
bundles  over spacetimes. A different $1$--cohomology has already proved 
useful in algebraic 
quantum field theory: a cohomology of the poset with values in a 
net of observables describes the 
superselection sectors. The formalism developed here can be adapted to 
this case.\smallskip

  We begin by explaining the notions of simplex, path and homotopy
in the context of posets showing that these notions behave in much
the same way as their better known topological counterparts.
We define the fundamental group of a path-connected poset which, in
practice, coincides with the fundamental group of the spacetime.
We then explain the $1$--cohomology of a poset with values in $G$
linking it to homotopy: the category of $1$--cocycles is equivalent
to the category of homomorphisms from the fundamental group to
$G$.\smallskip

  Having defined principal bundles, we next introduce the appropriate
notion of connection and curvature and investigate the set of
connections on a principal bundle, these being thus associated with
a particular $1$--cohomology. We discuss holonomy and prove a
version of the Ambrose-Singer Theorem. \smallskip

  We finally introduce the notion of gauge transformation and the
action of the group of gauge transformations on the set of connections
of a principal bundle. We also relate flat connections to homomorphisms 
from the fundamental group into $G$. We end by giving a brief outlook.

\section{Homotopy of posets}
\label{A}
We introduce some preliminary notions  
and results on posets.
We will start by defining the simplicial set associated with a poset  
and arrive at the notion of a simply connected poset. Throughout this 
section, we will consider a poset $K$ and denote its order relation by $\leq$. 
References for this section  are \cite{Rob3, Rob4, Ruz}.
\paragraph{The simplicial set of $K$ :} Underlying cohomology is what
is called the simplicial category $\Delta^+$  
that can be  realized  in various ways. The simplest way is to take the
objects of $\Delta^+$  to be the finite ordinals,
$\Delta_n=\{ 0,1,\ldots,n\}$   and  to  take  the  arrows to be 
the monotone mappings.  All these monotone mappings are compositions of two particular simple types of 
mapping; the injective monotone mappings from one ordinal to the succeeding 
ordinal denoted $d^n_i:(n-1)\to n$, with $i\in\{0,1,\ldots,n\}$,  
and defined as 
\[
 d^n_i(k) \defi \left\{
 \begin{array}{cl}
  k &    k<i\ ,  \\ 
  k+1 & otherwise\ ;
 \end{array}
 \right.
\]
and the 
surjective monotone mappings from one ordinal to the preceding one
denoted $s^n_i:(n+1)\to n$, with $i\in\{0,1,\ldots,n\}$, and defined as 
\[
  s^n_i(k) \defi \left\{
 \begin{array}{cl}
   k &   k\leq i\ ,  \\ 
  k-1 & otherwise\ .
 \end{array}
 \right.
\]
The superscripts of the symbols $d^n_i$ and $s^n_i$ are usually
omitted. The following identities allow one to compute effectively
$$d_id_j=d_{j+1}d_i,\,\, i\leq j;\quad s_js_i=s_is_{j+1},\,\,i\leq j;$$
$$s_jd_i=d_is_{j-1},\,\, i<j;\quad s_jd_j=s_jd_{j+1}=1;\quad 
  s_jd_i=d_{i-1}s_j,\,\, i>j+1.$$
Actually, each monotone map can be factorized uniquely as the composition of a 
surjective monotone map and an injective monotone map.\smallskip

We may also regard $\Delta_n$ as a partially ordered set, namely as 
the set of its non-void subsets ordered under inclusion. 
We denote this poset by $\widetilde{\Delta}_n$.
Any map, in particular a monotone one,  $m:\Delta_n\to\Delta_p$ induces, in an
obvious way, an order-preserving map of the partially ordered sets 
$\widetilde{\Delta}_n$ and $\widetilde{\Delta}_p$,  denoted by $\widetilde{m}$.
We can then define a \emph{singular $n$--simplex} of a poset $K$ 
to be an order preserving map $f:\widetilde{\Delta}_n\rightarrow K$. 
We denote the set of singular $n-$simplices  by $\Si_n(K)$, and call 
the \emph{simplicial set} of $K$ 
the set $\Si_*(K)$ of all singular simplices. Note that a map 
$m:\Delta_n\to\Delta_p$ induces 
a map $m^*:\Sigma_p(K)\to\Sigma_n(K)$, 
where $m^*(f)\defi f\circ \widetilde{m}$ with 
$f\in\Sigma_p(K)$. In particular, we have maps
\begin{equation}
\label{A:0}
\begin{array}{l}
\partial_i:\Sigma_n(K)\rightarrow \Sigma_{n-1}(K), \mbox{ where } \   
\partial_i\defi d_i^*, \\[8pt] 
\sigma_i :\Sigma_n(K)\rightarrow\Sigma_{n+1}(K),  \mbox{ where } \    
\sigma_i\defi s_i^*,
\end{array} 
\end{equation}
called \emph{boundaries} and \emph{degeneracies}, respectively. 
One can easily check the following relations
\begin{equation}
\label{A:1}
\begin{array}{c}
\partial_i  \partial_j = 
\partial_j  \partial_{i+1}, \ \  i\geq j;  \qquad 
\sigma_i\sigma_j=\sigma_{j+1}\sigma_i,  \ \  i\leq j;  \\[8pt]
\partial_i\sigma_j=\sigma_{j-1}\partial_i, \   i<j;
\ \   \partial_j\sigma_j=
\partial_{j+1}\sigma_j=1; \ \   \partial_i\sigma_j=\sigma_j\partial_{i-1},
\   i>j+1
\end{array}
\end{equation}
From now on, we will denote: the composition 
$\partial_i\partial_j$ by  the symbol $\partial_{ij}$;
$0$--simplices by the letter $a$; $1$--simplices by $b$; 
$2$--simplices by $c$ and a generic $n$--simplex by $d$. 
A $0$--simplex $a$ is just an element of the poset. 
Inductively, for $n\geq 1$, an $n-$simplex $d$ is formed
by $n+1$  $(n-1)-$simplices $\partial_0d, \ldots,\partial_{n}d$,
whose boundaries are constrained by the 
relations (\ref{A:1}), and by a $0$--simplex $|d|$ called the \emph{support}
of $d$ such that $|\partial_0d|, \ldots,|\partial_{n}d|\leq |d|$.
The ordered set $(\partial_0d,\ldots,\partial_nd)$, 
denoted $\partial d$, is called the \emph{boundary} of $d$. 
We say that an $n$--simplex is \emph{degenerate} if it is of the form
$\si_i(d)$ for some $(n-1)$--simplex $d$ and 
for some $i\in\{0,1,\ldots,n-1\}$. For instance, 
by using the relations (\ref{A:1}), 
it is easily seen that 
\[
 \partial_0\si_0(a)=\partial_1\si_0(a)= a, \ \ \ 
 |\si_0(a)|= a, 
\]
for any $0$--simplex $a$. In general, we have $|\si_i(d)|=|d|$.\smallskip

In later applications the following class of simplices will be
important. A $1-$simplex $b$ is said to be \emph{inflating} 
whenever 
\begin{equation}
\label{A:1a}
\partial_1b\subseteq \partial_0b.
\end{equation}
By induction for $n\geq 1$: an $n-$simplex $d$ is said to be 
\emph{inflating} whenever all its $(n-1)$--boundaries
$\partial_0d,\ldots,\partial_nd$ are inflating $(n-1)$--simplices.
For instance, if $c$ is  an  inflating $2$--simplex, then 
\[
\partial_{11}c=\partial_{12}c\subseteq 
\partial_{02}c= \partial_{10}c 
\subseteq \partial_{00}c=\partial_{01}c.
\] 
Any $0$-simplex will be regarded as inflating.
We will denote the set of inflating $n$--simplices by $\Si^{\inf}_n(K)$. 
Given a monotone mapping $m:\Delta_p\to\Delta_n$ then 
$m^*(d)=d\circ \widetilde{m}$ is 
inflating if $d$ is. Thus $\Sigma^{\inf}_*(K)$ is a simplicial subset of 
$\Sigma_*(K)$. \smallskip

We now deal with the notion of orientation of singular simplices 
of a poset $K$.  In general,  one says that a pair of simplices 
have  the same orientation whenever one can be obtained from  the 
other  by means of a even permutation of its vertices. 
The resulting equivalence relation will be called \emph{oriented equivalence}.
Notice that the $k^{th}$-vertex associated with 
an $n$--simplex $d$ is the $0$--simplex given by $\partial_{012\cdots
\widehat{k}\cdots n}d$, where $\widehat{k}$ means that the index 
$k$ is omitted. Given a permutation $\si$ of $(n+1)$ elements 
we denote by $d^\si$ the $n$--simplex obtained by permuting 
the vertices of $d$ according to $\si$ and leaving fixed the 
related supports. To be precise,
we define $d^\si$ as the $n$--simplex such that  $|d^\si|=|d|$
and  
\begin{equation}
\label{A:2}
  \partial_{012\cdots \widehat{k}\cdots n}d^{\si} = 
  \partial_{012\cdots \widehat{\si(k)}\cdots n}d, \qquad k\in
  \{0,1,\ldots, n\}
\end{equation}
These $n+1$ relations and the commutation relations (\ref{A:1}) allow 
one to compute how  the boundaries of $d^\si$ are related to those 
of $d$.  As an example, let  $\si$ be the transposition 
$(01)$. Then $a^{\sst{(01)}}=a$ for any $0$--simplex $a$. Inductively
for $n\geq 1$,  if $d$ is an $n$--simplex, then 
$|d^{\sst(01)}|= |d|$ and  
\begin{equation}
\label{A:3}
 \partial_0 d^{\sst{(01)}} =  \partial_0 d, \ \ 
 \partial_1 d^{\sst{(01)}} =  \partial_0 d, \ \ \  
 \partial_i d^{\sst{(01)}} =  
 (\partial_i d)^{\sst{(01)}} \ \ i\in \{2,3,\ldots,n\}. 
\end{equation}
Now, observe that the mapping 
\[
  \mathbb{P}(n+1)\times \Si_n(K) \ni (\si,d)\to d^\si\in\Si_n(K),
\]
where $\mathbb{P}(n+1)$ is the group of the permutations of $(n+1)$ elements, 
defines an action of $\mathbb{P}(n+1)$ on  $\Si_n(K)$. 
Two $n$--simplices  $d$ and $d_1$ are said to have 
the \emph{same orientation}  if  there exists an \emph{even} 
permutation $\si$ of $\mathbb{P}(n+1)$ such that $d_1 = d^\si$;
they have a \emph{reverse orientation} if there is 
an  \emph{odd} permutation $\si$  of $\mathbb{P}(n+1)$ such that 
$d_1 = d^\si$.   We denote 
by $[d]$ the equivalence class of $1$--simplices which have the same 
orientation as $d$, and  by $\overline{[d]}$ 
the equivalence class of $1$--simplices whose  orientation 
is the reverse of $d$.  
Notice that for any $0$--simplex $a$ we have $[a]=\overline{[a]}=a$. For
$1$--simplices we have $[b] = \{b\}$, while  
$\overline{[b]} =\{\overline{b}\}$,
where $\overline{b}$ is the $1$--simplex defined as 
\begin{equation}
\label{A:4}
 |\overline{b}| =|b|, \ \ \partial_0\overline{b}= \partial_1b, \ \ 
  \partial_1\overline{b}= \partial_0b. 
\end{equation}
In the following 
we will refer to $\overline{b}$ as the \emph{reverse} of the
$1$--simplex $b$ (note that $\overline{b}=b^{\sst{(01)}}$). 
For a $2$--simplex $c$ we
have 
\begin{equation}
\label{A:4a}
[c]=\{c, c^{\sst{(02)(01)}}, c^{\sst{(12)(01)}}\}, \ \ 
\overline{[c]}=\{c^{\sst{(01)}}, c^{\sst{(02)}}, c^{\sst{(12)}}\}.
\end{equation}
For instance,   $|c^{\sst{(02)(01)}}|=|c|$ and  
\begin{equation}
\label{A:5}
  \partial_0c^{\sst{(02)(01)}} =  \partial_2c,  \ \ 
  \partial_1c^{\sst{(02)(01)}} =  \overline{\partial_0c},  \ \ 
  \partial_2c^{\sst{(02)(01)}} =  \overline{\partial_1c};
\end{equation}
while $|c^{\sst{(12)(01)}}|=|c|$ and
\begin{equation}
\label{A:6}
  \partial_0c^{\sst{(12)(01)}} =  \overline{\partial_1c},  \ \ 
  \partial_1c^{\sst{(12)(01)}} =  \overline{\partial_2c},  \ \ 
  \partial_2c^{\sst{(12)(01)}} =  \partial_0c.
\end{equation}
In contrast to  the usual cohomological theories,
we do not  identify an $n$--simplex $d$ with its equivalence class 
$[d]$. This is  because in the following we will  deal with  
the curvature of a connection which is, in general,  not invariant 
under oriented equivalence.
\paragraph{Paths :} Given $a_0,a_1\in\Si_0(K)$, 
a  \emph{path from $a_0$ to $a_1$} is a finite ordered 
sequence $p=\{b_n,\ldots,b_1\}$ of $1$--simplices satisfying the relations
\[
 \partial_1b_1=a_0, \ \ \ \partial_0 b_{i} = \partial_1 b_{i+1} \  
            \mbox{ for } 
\ i\in\{1,\ldots,n-1\}, \ \ \  \partial_0b_n=a_1.
\]
The \emph{starting point} of $p$, written $\partial_1p$, 
is the $0$--simplex $a_0$, while the \emph{endpoint} of $p$, 
written $\partial_0p$, is the $0$--simplex  $a_1$. The \emph{boundary} 
of $p$ is the ordered set $\partial p\defi \{\partial_0p,\partial_1p\}$. 
A path $p$ is said to 
be a \emph{loop} if $\partial_0p=\partial_1p$. The \emph{support} $|p|$
of the path $p$ is the set 
\[
|p|\defi \{|b_1|,\ldots,|b_n|\}.
\]
We will denote the set of paths from $a_0$ to $a_1$ by $K(a_0,a_1)$, 
and the loops having endpoint $a_0$ by $K(a_0)$. $K$ will be 
assumed to be \emph{pathwise connected},  i.e. $K(a_0,a_1)$ is 
never void.   The set of paths is equipped with the following 
operations. Consider a path $p=\{b_n,\ldots,b_1\}\in K(a_0,a_1)$. 
The \emph{reverse} $\overline{p}$ is the path 
\[
\overline{p}\defi \{ \overline{b}_1,\ldots,\overline{b}_n\}\in K(a_1,a_0).
\]
The \emph{composition} of $p$ with a path 
$q=\{b'_k,\ldots b'_1\}$ of $K(a_1,a_2)$,  
is  defined by 
\[
q*p\defi \{b'_k,\ldots, b'_1,b_n,\ldots,b_1 \}\in K(a_0,a_2).
\]
Note that the reverse $-$ is involutive and 
the composition  $*$ is associative. In particular note 
that any path $p=\{b_n,\ldots,b_1\}$ can be also seen as the
composition
of its $1$--simplices, i.e., $p=b_n*\cdots*b_1$.  \\
\indent An \emph{elementary deformation} of a path $p$
consists in replacing a $1$--simplex $\partial_1c$ of the path
by a pair $\partial_0c,\partial_2c$, where $c\in\Si_2(K)$, or, conversely 
in replacing a consecutive pair $\partial_0c,\partial_2c$ of $1$--simplices 
of $p$ by a single $1$--simplex $\partial_1c$. Two paths with the same endpoints 
are \emph{homotopic} if they can be obtained from one other 
by a finite set of elementary deformations. Homotopy defines 
an equivalence relation  $\sim$  
on the set of paths with the same endpoints, which is compatible 
with reverse and composition,
namely 
\begin{equation}
\label{A:8}
\begin{array}{rclr}
 p\sim q   &   \iff   &    \overline{p}\sim\overline{q}, & 
    p,q\in  K(a_0,a_1); \\
 p\sim q, \ p_1\sim q_1  &   \Rightarrow  & 
  p_1*p\sim q_1* q, \ \   &  p_1,q_1\in  K(a_1,a_2).
\end{array}
\end{equation}
Furthermore, for any $p\in K(a_0,a_1)$, the following relations hold:
\begin{equation}
\label{A:9}
\begin{array}{rcl} 
  p* \sigma_0(a_0) \sim p & \mbox{ and } &  p \sim \sigma_0(a_1)*p;  \\
 \overline{p}* p\sim \si_0(a_0) & \mbox{ and }  & 
 \si_0(a_1)\sim p*\overline{p},
\end{array}
\end{equation}
where $\si_0(a_0)$ is the $1$--simplex degenerate to $a_0$. 
\paragraph{The first homotopy group: } Fix $a_0\in\Si_0(K)$, and define 
\[
 \pi_1(K,a_0)\defi K(a_0) /\sim,
\]
the quotient of the set of loops with endpoints $a_0$ 
by the homotopy equivalence relation. Let 
$[p]$ be the equivalence class associated with the loop $p\in K(a_0)$,
and let  
\[
 [p]\cdot [q] = [p*q], \qquad [p],[q]\in  \pi_1(K,a_0).
\]
$\pi_1(K,a_0)$ with this composition rule is a group: the identity 
is the equivalence class $[\sigma_0(a_0)]$ associated with the degenerate 
$1$--simplex $\sigma_0(a_0)$; the inverse of $[p]$ is the equivalence class 
$[\overline{p}]$ associated with the reverse $\overline{p}$
of $p$. $\pi_1(K,a_0)$ is the \emph{first homotopy group} 
of $K$ based on $a_0$. Since $K$ is pathwise connected 
$\pi_1(K,a_0)$ is isomorphic to $\pi_1(K,a)$ for any 
$a\in\Si_0(K)$; this isomorphism class is the
\emph{fundamental group} of $K$, written $\pi_1(K)$. 
If $\pi_1(K)$ is trivial, then  $K$ is said to 
be \emph{simply connected}. It turns out that 
if $K$ is directed\footnote{The poset $K$ is directed whenever for any pair
$\dc,\dc_1\in K$ there exists $\dc_2\in K$ such that $\dc,\dc_1\leq
\dc_2$}, then $K$ is simply connected.\\
\indent The link between the homotopy group of a poset and 
the corresponding topological notion,
can be achieved as follows. Let $M$ be an arcwise connected 
manifold and let 
$K$ be a base for the topology of $M$ whose elements
are arcwise  and simply connected, open  subsets of $M$.
Consider the poset formed by ordering $K$ under \emph{inclusion}. 
Then $\pi_1(M)=\pi_1(K)$,  
where $\pi_1(M)$ is the fundamental group of $M$. 
\paragraph{Coverings :} A partially ordered set 
$K$ can be equipped with a $\mathrm{T}_0$  topology called 
the Alexandroff topology.   
In this topology, a subset $U\subseteq K$ is said to be 
\emph{open} whenever given $\dc\in U$ and $\dc_1\in K$,  
if $\dc\leq\dc_1$  then  $\dc_1\in U$.  An  \emph{open covering}  of $K$ 
is a family $\Uc$ of open sets $U$ of $K$ 
such that for any  $\dc\in K$ there is $U\in\Uc$ with $\dc\in U$.
A particular  covering is that formed by the collection 
$\{U_a , \  a\in\Si_0(K)\}$  of open sets of $K$ defined by
\begin{equation}
\label{A:10}
  U_a \defi \{\dc\in K \ | \ a\leq \dc\}, \qquad  a\in\Si_0(K).
\end{equation}
We call this covering the \emph{fundamental  covering} of $K$.
Note that if $\Uc$ is an open  covering of $K$, then for any $0$--simplex $a$
there is $U\in\Uc$ such that $U_a\subseteq U$.
\section{Cohomology of posets}             
\label{B}                                  
The present section deals with the, in general non-Abelian,  
cohomology of a pathwise connected 
poset $K$ with values in a group $G$.
The first part is devoted to explaining the motivation
for studying the non-Abelian cohomology of a poset and 
to  defining  an $n$--category. The general theory 
is developed in the  second part: we introduce 
the set of $n$--cochains, for $n=0,1,2,3$,   the coboundary 
operator, and the cocycle identities up to the $2^{nd}$-degree.   
In the last part we study the 
$1$--cohomology, in some detail, relating it to the first 
homotopy group of a poset. 
\subsection{Preliminaries}
\label{Ba}

The cohomology of the poset $K$ with values in an Abelian group $A$, 
written additively, 
is the cohomology of the set of singular simplices 
$\Si_*(K)$  with values in  $A$. To be precise,  
one can define  the set $\Crm^n(K,A)$ of $n$--cochains of $K$ 
with values in $A$ 
as the set of functions $v:\Si_n(K)\rightarrow A$. 
The coboundary operator $\dr$ defined by 
\[
\dr v (d) =  \sum^n_{k=0} (-1)^k  \ v(\partial_k d), \qquad d\in\Si_n(K),
\]
is a  mapping $\dr:\Crm^n(K,A)\rightarrow\Crm^{n+1}(K,A)$ satisfying
the equation $\dr\dr v = \io$, 
where $\io$ is the trivial cochain. This allows one
to define the $n$--cohomology groups. For a  non-Abelian 
group $G$ no choice of ordering gives the identity 
$\dr\dr v= \io$.\\
\indent One motivation for studying the cohomology of a poset
$K$ with values in a non-Abelian group comes from algebraic  quantum 
field theory. The leading idea of this approach is that all 
the physical content of a quantum system  is encoded
in the observable net, an inclusion preserving 
correspondence which associates to any 
open and bounded region of  Minkowski space
the algebra  generated by the observables measurable within that
region. The collection of these regions forms a poset when ordered 
under inclusion.  A $1$--cocycle equation  arises in  studying 
charged  sectors of the  observable  net: the charge transporters of 
sharply localized charges are  $1$--cocycles of the poset taking  values 
in the group of  unitary operators of the observable 
net \cite{Rob1}. 
The attempt to include more general charges in the framework 
of local quantum physics, charges of electromagnetic type in
particular, has led one  to derive higher cocycles 
equations, up to the third degree \cite{Rob2,Rob3}. The difference, 
with respect to the Abelian case, is that a $n$--cocycle equation 
needs $n$--composition laws. Thus in non-Abelian cohomology 
instead, for example, of trying to take coefficients in a non-Abelian
group the $n$--cocycles take values in an $n$--category associated
with the group. The cocycles equations can be understood 
as pasting together simplices, and, in fact, a $n$--cocycle 
can be seen as a representation in an $n$--category 
of the algebra of an  oriented  $n$--simplex \cite{Str}.\smallskip 
    
Before trying to learn the notion of an $n$--category, it helps to 
recall that a category can be defined in two equivalent manners. 
One definition is based on the set of objects and the corresponding
set of arrows. However, it is possible to define a category referring
only to the set of arrows. Namely, a category 
is a set $\C$, whose elements are called arrows,  having a partial 
and associative composition law $\diamond$, 
and such that any element of $\C$ has left and right $\diamond$-units. 
This amounts to saying that (i) $(f\diamond g)\diamond  h$ is defined 
if, and only if, $f\diamond(g\diamond h)$ is defined and they are equal;
(ii) the triple $f\diamond g\diamond h$ is defined if, and only if, $f\diamond g$ and 
     $g\diamond h$ are defined; 
(iii) any arrow $g$ has  a left and a right unit $u$ and $v$, 
   that is $u\diamond g= g$ and $g\diamond v = g$. 
In this formulation the set of objects are the set of units.\\  
\indent An $n$--category is a set $\C$ with an ordered set 
of $n$ partial composition laws. This means that 
$\C$ is a category with respect to any such 
composition law $\diamond$. Moreover,
if $\times$ and $\diamond$ are two such composition laws 
with $\times\prec \diamond$ then: 
\begin{enumerate}
\item every $\times$-unit is a $\diamond$-unit; 
\item $\times$-composition  of $\diamond$-units, when defined, 
      leads to $\diamond$-units;
\item the following relation, called the \emph{interchange law}, holds: 
\[
        (f\times h) \diamond  (f_1\times h_1) = 
       (f\diamond f_1) \times (h\diamond h_1),  
\]
whenever the right hand side is defined.
\end{enumerate}
An arrow $f$ is said to be a \emph{$k$--arrow}, for $k\leq n$,  if 
it is a unit for the $k+1$ composition law. 
To economize on brackets, from now on we adopt the convention 
that if $\times \prec\diamond$, then a $\times$--composition law 
is to be evaluated before a $\diamond$--composition. For example, 
the interchange law reads 
\[
        f\times h \diamond  f_1\times h_1 = 
       (f\diamond f_1) \times (h\diamond h_1).
\]
It is surprising that with this convention all the brackets
disappear from the coboundary equations (see below).\smallskip  

That an $n$--category is the right set of coefficients 
for a non-Abelian cohomology can be understood by the following
observation. Assume that $\times$ is Abelian, that is, 
$f\times g$ equals $g\times f$ whenever the compositions are defined.   
Assume that $\diamond$--units are $\times$--units. 
Let $1,1'$ be, respectively,  a left and a right $\diamond$--unit for
$f$ and $g$. 
By using the interchange law we have
\[
f\diamond g =  1\times f \diamond g\times 1'  =
 (1\diamond g)  \times  (f\diamond 1')  = g\times f.  
\]
Hence $\diamond$ equals $\times$ and both  composition laws are
Abelian. Furthermore, if $\star$ is a another composition law 
such that $\times \prec \star\prec \diamond$, then $\times=\star=\diamond$.

\subsection{Non-Abelian cohomology}
\label{Bb} 
The first aim is to introduce an $n$--category  
associated with a group $G$ to be  used 
as set of coefficients for the cohomology of the poset $K$.
To this end, we draw on a general procedure \cite{Rob6} 
associating to any $n$--category $\C$ where the  $n$--arrows are
invertible, with respect to any composition law, 
an $(n+1)$--category $\mathcal{I}(\C)$ with the same property.
This construction 
allows one to define 
the $(n+1)$--coboundary of a $n$--cochain in $\C$ as an
$(n+1)$--cochain in $\mathcal{I}(\C)$, at least for $n=0,1,2$.\\ 
\indent Before starting to describe non-Abelian cohomology, 
we introduce some  notation. The elements
of a group $G$ will be indicated  by  Latin letters.
The composition of two elements $g,h$ of $G$ will be denoted 
by $gh$, and by $e$ the identity of $G$.
Let $Inn(G)$ be the group of inner automorphisms
of $G$. We will use Greek letters to indicate the elements of
$Inn(G)$. By  $\alpha\tau$ we will denote  the inner automorphism of $G$
obtained by the composing  $\alpha$ with $\tau$, that is 
$\alpha\tau(h) \defi  \alpha(\tau(h))$ for any $h\in G$. The identity 
of this group, the trivial automorphism,  will be indicated by $\io$. 
Finally 
given $g\in G$, the equation 
\[
  g\, \al = \tau\, g
\]
means $g\alpha(h) = \tau(h)g$ for any $h\in G$. 
\paragraph{The categories $nG$ :} 
In degree $0$,  this is simply the group $G$ considered as a set. 
In  degree $1$ 
it is the category $1G$ having a single object, the group $G$, and as arrows 
the elements of the group. Composition of arrows is the composition 
in $G$. So we identify this category with $G$. 
Observe that the arrows of $1G$ are invertible. 
By applying the procedure provided in \cite{Rob6} 
we have that $\mathcal{I}(1G)$ is a $2$--category, denoted by 
$2G$, whose set of arrows is 
\begin{equation}
\label{Bb:1}
 2G\defi \{ (g,\tau) \ | \ g\in G, \ \ \tau\in Inn(G)\},  
\end{equation}
and whose  composition laws are defined by 
\begin{equation}
\label{Bb:2}
\begin{array}{rcll}
(g,\tau)\times  (h,\ga) & \defi & 
(g\tau(h),\tau\ga), \\
(g,\tau)\diamond  (h,\ga) & \defi & 
(gh,\ga), &  \mbox{ if }  \si_h\ga = \tau,
\end{array}
\end{equation}
where $\si_h$ is the inner automorphism associated with $h$. 
Some observations on $2G$ are in order.  
Note that the composition $\times$ is always defined. 
Furthermore, the set of $1$--arrows is the set of those elements 
of $2G$ of the form 
$(e,\tau)$. Finally, all the $2$--arrows are invertible. 
We can now  construct the $3$--category $\mathcal{I}(2G)$,
denoted by $3G$. It turns out that $3G$ is the set  
\begin{equation}
\label{Bb:3}
 3G\defi \{ (g,\tau,\ga) \ | \ g\in\mathcal{Z}(G), \ \ \tau,\ga\in
 Inn(G)\},  
\end{equation}  
where $\mathcal{Z}(G)$ is the centre of  $G$, with the following 
three composition laws
%
%
\begin{equation}
\label{Bb:4}
\begin{array}{rcll}
(g,\tau,\ga)\times  (g',\tau',\ga') & \defi & 
(gg',\tau\tau', \ga\tau\ga'\tau^{-1}), \\
(g,\tau,\ga)\diamond (g',\tau',\ga') & \defi & 
(gg', \tau', \ga\ga',), & \mbox{ if } \tau= \ga'\tau'\\
(g,\tau,\ga)\cdot  (g',\tau',\ga') & \defi & 
(gg', \tau,\ga), & \mbox{ if } \tau= \tau', \ \ga= \ga'.
\end{array}
\end{equation}
Note that $\cdot$ is Abelian.  The set 
of $1$-arrows $(3G)_1$ is the subset of  elements 
of $3G$ of the form 
$(e,\ga,\io)$, where $\io$ denotes the identity automorphism;
$2$--arrows $(3G)_2$ are  the elements of $3G$ of the form 
$(e,\tau,\ga)$. Finally, if 
$G$ is  Abelian, then $\times=\diamond=\cdot$ and the categories 
$2G$ and $3G$ are nothing but that the group $G$.
%
\paragraph{The set $\Si_n(K,G)$ of $n$--cochains :}
The next goal  is to define the set of $n$--cochains. Concerning 
$0$-- and $1$--cochains nothing change with respect to the Abelian case,
i.e.,  $0$--  and $1$--\emph{cochains} are, respectively 
functions $v:\Si_0(K)\to G$ and $u:\Si_1(K)\to G$. 
A \emph{$2$--cochain} $w$ is a pair of mappings $(w_1,w_2)$, where  
$w_i:\Si_i(K)\to (2G)_i$, for $i=1,2$ enjoying the relation
\begin{equation}
\label{Bb:5}
  w_2(c) \diamond  w_1(\partial_1c) =  
  w_1(\partial_0c) \times  w_1(\partial_2c) \diamond w_2(c), \qquad c\in\Si_2(K).
\end{equation}
This equation and the definition of the composition laws 
in $2G$  entail that a $2$--cochain $w$ is of the form 
\begin{equation}
\label{Bb:6}
\begin{array}{rcll}
 w_1(b)& = & (e, \tau_{b}), & b\in\Si_1(K), \\[3pt]
 w_2(c)& = & (v(c), \tau_{\partial_1c}), &  c\in\Si_2(K),
\end{array}
\end{equation}
where $v:\Si_2(K)\to G$, $\tau:\Si_1(K)\to Inn(G)$ are mapping 
satisfying the 
equation\footnote{Equation \ref{Bb:7} means 
that $v(c)$ intertwines 
$\tau_{\partial_1c}$ with $ \tau_{\partial_0c} \ \tau_{\partial_2c}$,
that is $v(c) \ \tau_{\partial_1c}(h) = \tau_{\partial_0c}(\tau_{\partial_2c}(h))
 \ v(c)$ for any, $h\in G$.}
\begin{equation}
\label{Bb:7}
 v(c) \ \tau_{\partial_1c} = \tau_{\partial_0c} \, \tau_{\partial_2c} \ v(c),
 \qquad c\in\Si_2(K). 
\end{equation}
This can be easily shown. In fact, according to the definition of 
$2G$ a $2$--cochain $w$ is of the form 
$w_1(b) = (e, \tau_b)$ for $b\in\Si_1(K)$, and  
$w_2(c) = (v(c), \be_c)$ for $c\in\Si_2(K)$. Now, 
the l.h.s. of  equation (\ref{Bb:5}) is defined if, and only if,  
$\tau_{\partial_1c} = \be_c$ for any $2$--simplex $c$.  
This fact and  equation (\ref{Bb:5})  entail  (\ref{Bb:7}) and 
(\ref{Bb:6}), completing the proof. \\
\indent A \emph{3-cochain} $x$ is $3$--tuple $(x_1,x_2,x_3)$ where 
$x_i:\Si_1(K)\to (3G)_i$, for $i=1,2,3$, satisfying the following 
equations 
\begin{equation}
\label{Bb:8a}
x_2(c) \diamond  x_1(\partial_1c)  =    x_1(\partial_0c) \times 
  x_1(\partial_2c) \diamond x_2(c),
\end{equation}
for any $2$--simplex $c$, and 
\begin{align}
\label{Bb:8b}
x_3(d) \cdot  x_1(\partial_{01}d)\times x_2&(\partial_3d) \diamond
  x_2(\partial_1d)    = \notag \\
   & =    x_2(\partial_0d) \times 
  x_1(\partial_{23}d)\diamond  x_2(\partial_2d) \cdot x_3(d), 
\end{align}
for any $3$--simplex $d$.
Proceeding as above, these equations and the composition laws of $3G$ 
entail  that a $3$--cochain $x$ is of the form 
\begin{equation}
\label{Bb:9}
\begin{array}{ll}
 x_1(b)  =(e,\ \tau_b,\ \io),  & b\in\Si_1(K),\\[3pt]
 x_2(c)  = (e,\  \tau_{\partial_1c}, \ \ga_c), 
     & c\in\Si_2(K),\\[3pt]
 x_3(d)  = (v(d), \ \tau_{\partial_{12}d},\  
           \gamma_{\partial_0d}\,  \gamma_{\partial_1d}), & d\in\Si_3(K),
\end{array}
\end{equation}
where $\tau:\Si_1(K)\to Inn(G)$, $v:\Si_3(K)\to Z(G)$, 
while $\ga:\Si_2(K)\to Inn(G)$ is the mapping defined as 
\begin{equation}
\label{Bb:9a}
 \ga_c \defi  \tau_{\partial_0c}\, \tau_{\partial_2c}\, \tau_{\partial_1c}^{-1}, 
 \qquad  c\in\Si_2(K).
\end{equation}
Note, in particular that $\ga_c \ \tau_{\partial_1c}
= \tau_{\partial_0c}\tau_{\partial_2c}$ for any $2$--simplex $c$.
This concludes the definition of the set of cochains. We will denote 
the set of $n$--cochains of $K$, for $n=0,1,2,3$,  
by $\Crm^n(K,G)$.\smallskip

Just a comment about the   definition of $1$--cochains. 
Unlike the usual cohomological theories  
$1$--cochains are neither required to be invariant under 
oriented equivalence of simplices 
nor to act trivially on degenerate simplices. However, as we
will see later, $1$--cocycles and connections fulfil these
properties.
\paragraph{The coboundary  and the cocycle identities:}
The next goal is to define the coboundary operator $\dr$.
Given a $0$--cochain $v$, then 
\begin{equation}
\label{Bb:10}
  \dr v(b)  \defi   v(\partial_0b) \ v(\partial_1b)^{-1}, 
\qquad b\in\Si_1(K).
\end{equation}
Given a $1$-cochain $u$,  then 
\begin{equation}
\label{Bb:11}
\begin{array}{ll}
  (\dr u)_1 (b)  \defi  (e, \ \mathrm{ad}(u(b))), & b\in\Si_1(K), \\[3pt]
  (\dr u)_2 (c)  \defi  (w_u(c), \ 
    \mathrm{ad}(u(\partial_1c))), \ \ &  c\in\Si_2(K),
\end{array}
\end{equation}
where $w_u$ is the mapping from $\Si_2(K)$ to $G$ defined as 
\begin{equation}
\label{Bb:11a}
w_u(c) \defi  u(\partial_0 c) \  
                           u(\partial_2 c) \  u(\partial_1c)^{-1},
                           \qquad c\in\Si_2(K).
\end{equation}
Finally, given a $2$-cochain $w$ of the form (\ref{Bb:6}),  then 
\begin{equation}
\label{Bb:12}
\begin{array}{ll}
(\dr w)_1 (b)  \defi  \big(e, \ \tau_b,\ \io\big),  &   b\in\Si_1(K), \\[3pt]
(\dr w)_2 (c)  \defi  \big(e, \  
   \tau_{\partial_1c}, \ \ga_c \big), &  c\in\Si_2(K), \\[3pt] 
  (\dr w)_3 (d)  \defi  \big( x_w (d), \ 
   \tau_{\partial_{12}d}, \ 
             \ga_{\partial_0d}\, \ga_{\partial_2d}\big), \ \ & 
 d\in\Si_3(K),
\end{array}
\end{equation}
where $\ga$ is the function from $\Si_2(K)$ to $Inn(G)$ defined 
by $\tau$ as in (\ref{Bb:9a}), and $x_w$ is the mapping 
$x_w: \Si_3(K)\to \mathcal{Z}(G)$ defined as 
\begin{equation}
\label{Bb:12a}
x_w(d)\defi v(\partial_0d) \, v(\partial_2d) \,
   \big(\tau_{\partial_{01}d}(v(\partial_3d))\,  v(\partial_1d)\big)^{-1}
\end{equation}
for any $3$--simplex $d$. Now, we call the \emph{coboundary operator} 
$\dr$ the mapping 
$\dr:\Crm^n(K,G)\to \Crm^{n+1}(K,G)$ defined 
for $n=0,1,2$ by the equations (\ref{Bb:10}),  (\ref{Bb:11}) and  
(\ref{Bb:12}) respectively.  This definition is well 
posed as shown by the following
\begin{lemma} 
\label{Bb:13}
For $n=0,1,2$,  the coboundary operator $\dr$ is a mapping 
$\dr:\Crm^n(K,G)\to \Crm^{n+1}(K,G)$, such that 
\[ 
 \dr\dr v\in ((2+k)G)_{k+1}, \qquad  v\in\Crm^k(K,G)
\] 
for $k=0,1$.
\end{lemma}
\begin{proof}
The proof of the first part of the statement
follows easily from the definition of $\dr$, except  that
the function $x_w$, as defined in (\ref{Bb:12a}),  
takes values in $\mathcal{Z}(G)$.
Writing,  for brevity,  $v_{i}$ for  $v(\partial_{i}d)$ and 
$\tau_{ij}$ for $\tau_{\partial_{ij}}$, and  using relations 
(\ref{A:1}) and  equation (\ref{Bb:7}) we have 
\[
v_0  \, v_2  \, \tau_{12}  =   v_0  \, \tau_{02}\, \tau_{22} \, v_2
 = v_0 \, \tau_{10}\, \tau_{22} 
  \, v_2
 =  \tau_{00} \, \tau_{20} \,  
\tau_{23} \,  v_0  \,   v_2,
\]
moreover
\begin{align*}
\tau_{01}(v_3) \, v_1  \,  
 \tau_{12} & = 
  \tau_{01} (v_3) \, 
         v_1 \, \tau_{11} 
 = \tau_{01}(v_3d) \, 
     \tau_{01} \, \tau_{21} \,  
         v_1 \\
&   = \tau_{01}(v_3d\, \tau_{13}) \, 
          v_1 
 = \tau_{01}(\tau_{03} \, \tau_{23}\, 
          v_3) \, 
          v_1 \\
&   = \tau_{01}\, \tau_{03}\,
 \tau_{23}\, \tau_{01}(v_3) \, v_1,
\end{align*}
Hence both $v(\partial_0d) \, v(\partial_2d)$ and 
$ \tau_{\partial_{01}d}(v(\partial_3d))\, v(\partial_1d)$ intertwine 
from $\tau_{\partial_{12}d}$ to $\tau_{\partial_{01}d}\, 
\tau_{\partial_{03}d}\, \tau_{\partial_{23}d}$. This entails that they 
differ  only by an element of $\mathcal{Z}(G)$, and this proves 
that $x_w$ takes values in $\mathcal{Z}(G)$. 
Now, it is very easy to see that 
$\dr\dr v\in (2G)_1$,  for any $0$--cochain $v$. 
So, let us prove that $\dr\dr u\in (3G)_2$ 
for any $1$--cochain $u$. Note that 
\[
\begin{array}{ll}
(\dr u)_1 (b) = (e,\mathrm{ad}(u(b))),  &  b\in\Si_1(K) \\[3pt]
(\dr u)_2 (c) = (w_u(c),
                      \mathrm{ad}(u(\partial_1c))),
                       &  c\in\Si_2(K),
\end{array}
\]
where $w_u$ is defined by (\ref{Bb:11a}). Then the proof follows 
once we have shown that 
\[
 w_u(\partial_0d) \ w_u(\partial_2d)  = 
 \mathrm{ad}(u(\partial_{01}d))\big(w_u(\partial_3d)\big) \ 
                           w_u(\partial_1d), \qquad (*) 
\]
for any $3$--simplex $d$. In fact, by  (\ref{Bb:12}) this identity 
entails that 
\[
\begin{array}{ll}
(\dr\dr u)_1 (b) = (e,\mathrm{ad}(u(b)), \io),  &  b\in\Si_1(K) \\[3pt]
(\dr\dr u)_2 (c) = (e,\ \mathrm{ad}(u(\partial_1c)), \
                       \mathrm{ad}(w_u(c)))
                       &  c\in\Si_2(K),\\[3pt]
(\dr\dr u)_3 (d) = (e,\ \mathrm{ad}(u(\partial_{12}d)), \ 
      \mathrm{ad}\big(w_u(\partial_0d) w_u(\partial_2d)\big)), &  
 d\in\Si_{3}(K). 
\end{array}
\]
So let us prove $(*)$. 
Given $d\in\Si_3(K)$ and using  relations (\ref{A:1}), we have 
\begin{align*}
\mathrm{ad}(u(\partial_{01}d))\big(w_u(\partial_3d)\big) \ & 
w_u(\partial_1d)  = \\
& =  u(\partial_{01}d) \,  w_u(\partial_3d) \, u(\partial_{01}d)^{-1} \, 
    w_u(\partial_1d)  \\
& =   u_{01}\, ( u_{03}\,  u_{23}\,  u_{13}^{-1})\,  u_{01}^{-1}\, 
       ( u_{01}\,  u_{21}\,  u_{11}^{-1}) \\
& =   u_{01}\, u_{03}\, u_{23}\,  u_{11}^{-1} \\
& =   u_{01}\, u_{03}\, u_{02}^{-1}\, u_{02}\, u_{23}\, u_{11}^{-1} \\
& =   u_{00}\, u_{20}\, u_{10}^{-1}\, u_{02}\, u_{22}\, u_{12}^{-1} \\
& =  w_u(\partial_0d)\, w_u(\partial_2d),
\end{align*}
where we have used the  notation introduced above. 
This completes  the proof.
\end{proof}
In words this lemma says that if $v$ is a $0$--cochain, then 
$\dr\dr v$ is  a  $2$--unit  of $2G$;
if $u$ is a $1$--cochain, then 
$\dr\dr v$ is a  $3$--unit  of  $3G$.\smallskip

We now are in a position to introduce 
the definition of  an $n$--cocycle. 
\begin{df}
\label{Bb:14}
For $n=0,1,2$,  an $n$--cochain $v$ is said to be an 
\textbf{$n$--cocycle} whenever 
\[
 \dr v\in \big((n+1)G\big)_n.
\]
It is said to be an \textbf{$n$--coboundary} whenever 
\[
 v\in \dr\big((n-1)G\big)
\]
(for $n=0$ this means that $v(a) = e$ for any
$0$--simplex $a$). 
We will denote the set of $n$--cocycles by $\Zrm^n(K,G)$,
and by $\Brm^n(K,G)$ the set of $n$--coboundaries.
\end{df} 
Lemma \ref{Bb:13} entails  that
$\Brm^n(K,G)\subseteq \Zrm^n(K,G)$ 
for $n=0,1,2$. Although it is outside  the scope  of this paper, 
we note that this relation holds also for $n=3$. One can check 
this assertion by using the $3$--cocycle given in \cite{Rob2}.

It is very easy now to derive the cocycle equations. 
A $0$--cochain $v$ is a \emph{$0$--cocycle} if
\begin{equation}
\label{Bb:15}
  v(\partial_0b) = v(\partial_1b), \qquad b\in\Si_1(K).
\end{equation}
A $1$--cochain $z$ is a \emph{$1$--cocycle} if 
\begin{equation}
\label{Bb:16}
  z(\partial_0c) \, z(\partial_2c) = z(\partial_1c), \qquad c\in\Si_2(K).
\end{equation}
Let  $w = (w_1,w_2)$ be $2$--cochain of the form 
$w_1(b) =  (e, \tau_{b})$ for $b\in\Si_1(K)$,
$w_2(c) =  (v(c), \tau_{\partial_1c})$ for $c\in\Si_2(K)$, 
where $v$ and $\tau$ are mappings satisfying (\ref{Bb:7}). 
Then $w$ is a \emph{$2$--cocycle} if 
\begin{equation}
\label{Bb:17}
 v(\partial_0d)\,  v(\partial_2d)  = 
 \tau_{\partial_{01}d}\big(v(\partial_3d)\big)\, 
                          v(\partial_1d), \qquad d\in\Si_3(K).
\end{equation}   
In the following we will mainly deal with $1$--cohomology. 
Our purpose will be to show 
that the notion of $1$--cocycle admits an interpretation 
as a principal bundle over a poset and that this kind of bundle 
admits connections. The $2$--coboundaries 
enter the game as the curvature of connections. 
Since the poset is pathwise connected, it turns out that 
any $0$--cocycle $v$ is a constant function.  
Thus the $0$--cohomology of $K$ yields no 
useful information. 
\subsection{$1$--Cohomology}
\label{Bc}
This section is concerned with  $1$--cocycles of the poset $K$.
In the first part we introduce some basic notions 
that will be used throughout this paper. The second part deals 
with $1$--cocycles. We will derive some results confirming the 
interpretation of a $1$--cocycle as a principal bundle 
over a poset.
This interpretation will become clear in Section \ref{C}
In the last part we discuss the connection between 
$1$--cohomology and homotopy of posets.     
\paragraph{The category of $1$-cochains :} Given a $1$-cochain 
$v\in\Crm^1(K,G)$, we can and will extend $v$ 
from $1$-simplices to paths by defining for $p=\{b_n,\ldots,b_1\}$ 
\begin{equation}
\label{Bc:1}
v(p) \defi v(b_n)\, \cdots\,  v(b_2)\, v(b_1).
\end{equation}
\begin{df}
\label{Bc:2}
Consider $v,v_1\in\Crm^1(K,G)$.  
A \textbf{morphism} $\fr$ from $v_1$ to $v$ 
is a function $\fr:\Si_0(K)\rightarrow G$  satisfying the equation 
\[
  \fr_{\partial_0p} \,  v_1(p)  = v(p) \, \fr_{\partial_1p}, 
\]
for all paths $p$. We denote the set of morphisms from 
$v_1$ to $v$ by $(v_1,v)$. 
\end{df} 
There is an obvious composition law between morphisms given by pointwise 
multiplication and this makes $\Crm^1(K,G)$ into a category. 
The identity arrow $1_v\in(v,v)$ takes the constant value 
$e$, the identity of the group. 
Given a group homomorphism $\gamma:G_1\to G$ and a 
morphism $\fr\in(v_1,v)$ of 
$1$--cochains with values in $G_1$ then $\gamma\circ v$, defined as  
\begin{equation}
\label{Bc:3a}
(\gamma\circ  v) (b)\defi \ga(v(b)), \qquad b\in\Si_1(K), 
\end{equation}
is a $1$--cochain with values in $G$, and $\gamma\circ\fr$ defined as 
\begin{equation}
\label{Bc:3b}
(\gamma\circ\fr)_a\defi \ga(\fr_a), \qquad a\in\Si_0(K), 
\end{equation}
is a morphism of $(\gamma\circ v_1,\gamma\circ v)$. 
One checks at once that $\gamma\circ$ is a functor from $\Crm^1(K,G_1)$ to 
$\Crm^1(K,G)$, and that if $\ga$ is a group isomorphism, then 
$\ga\circ$ is an isomorphism of categories.\\
\indent Note that $\fr\in(v_1,v)$ implies $\fr^{-1}\in(v,v_1)$, where
$\fr^{-1}$ here denotes the composition of $\fr$ with the inverse of 
$G$. We say that $v_1$ and $v$ are \emph{equivalent}, 
written $v_1\cong v$, whenever $(v_1,v)$ is nonempty. 
Observe that  a $1$--cochain $v$ is equivalent 
to the trivial $1$--cochain $\imath$ if, and only if, it is 
a $1$--coboundary.
We will say that $v\in \Crm^1(K,G)$  is 
\emph{reducible} if there exists a proper subgroup 
$G_1\subset G$ and a $1$--cochain $v_1\in\Crm^1(K,G_1)$ with  
$\gamma\circ v_1$ equivalent to $v$, where $\gamma$ denotes the 
inclusion $G_1\subset G$. If $v$ is not reducible it will be said 
to be \emph{irreducible}.  
 
A $1$--cochain $v$ is said to be \emph{path-independent} 
whenever given a pair of paths $p,q$, then 
\begin{equation}
\label{Bc:4}
 \partial p =  \partial q \ \ \Rightarrow \ \ v(p)=v(q) \ .
\end{equation} 
Of course, if $v$ is path-independent then so is any equivalent 
$1$--cochain.
It is worth observing that if $\gamma$ is an injective homomorphism 
then $v$ is path-independent if, and only if, $\gamma\circ v$ is 
path-independent. 
\begin{lemma}
\label{Bc:5}
Any $1$--cochain is path-independent if, and only if, it is a $1$--coboundary.
\end{lemma}
\begin{proof}
Assume that $v\in\Crm^1(K,G)$ is path-independent.  Fix a $0$--simplex 
$a_0$. For any $0$--simplex $a$, choose a path $p_a$ from $a_0$ to $a$
and define $\fr_a\defi v(p_a)$.
As $v$ is path-independent, for any $1$--simplex $b$ we have
\[
v(b)\, \f_{\partial_1b}= v(b) \,  v(p_{\partial_1b})=
v(b*p_{\partial_1b})=v(p_{\partial_0b})=\f_{\partial_0b}.
\]
Hence $v$ is a $1$--coboundary, see \ref{Bb:10}. 
The converse is obvious.
\end{proof}
\paragraph{$1$-Cocycles as principal bundles :} Recall that 
a $1$--cocycle $z\in\Zrm^1(K,G)$ is a mapping 
$z:\Si_1(K)\to G$ satisfying the equation 
\[
 z(\partial_0c)\, z(\partial_2c) = z(\partial_1c), \qquad 
c\in\Si_2(K)
\]
Some  observations are in order. 
First, the trivial $1$--cochain $\imath$ is a 
$1$--cocycle (see Section
\ref{Ba}).  So, from now on, we will refer to $\imath$ as the \emph{trivial 
$1$-cocycle}. Secondly,  if $z$ is a $1$--cocycle  then so is any
equivalent $1$--cochain. In fact, let  
$v\in\Crm^1(K,G)$ and let $\fr\in (v,z)$. 
Given a $2$--simplex $c$ we have 
\begin{align*}
v(\partial_0c)\,  v(\partial_2c) & = 
\fr^{-1}_{\partial_{00}c}\,  z(\partial_0c)\,  
\fr_{\partial_{10}c}\,  
\fr^{-1}_{\partial_{02}c}\,  z(\partial_2c)\,  
\fr_{\partial_{12}c} \\
& = \fr^{-1}_{\partial_{00}c}\,  z(\partial_0c) \,  z(\partial_2c)\,  
\fr_{\partial_{12}c} 
= \fr^{-1}_{\partial_{00}c}\,  z(\partial_1c)\, 
\fr_{\partial_{12}c} \\ 
& = 
\fr^{-1}_{\partial_{01}c}\,  z(\partial_1c)\, 
\fr_{\partial_{11}c} = v(\partial_1c),
\end{align*}
where relations (\ref{A:1}) have been used.
\begin{lemma}
\label{Bc:6}
Let $\gamma:G_1\to G$ be a group homomorphism. Given 
$v\in \Crm^1(K,G_1)$ consider $\ga\circ v\in\Crm^1(K,G)$. 
Then: if $v$ is a $1$--cocycle, 
then $\ga\circ v$ is a $1$--cocycle; the converse holds 
if $\ga$ is injective.
\end{lemma}
\begin{proof} 
If $v$ is a $1$--cocycle, it is easy to see that 
$\ga\circ v$ is a $1$--cocycle too. Conversely,  assume 
that $\ga$ is injective and that $\ga\circ v$ is a $1$--cocycle, then 
\[
\ga\big(v(\partial_0c)\  v(\partial_2c)\big)  = 
\ga\circ v(\partial_0c)\  \ga\circ v(\partial_2c) =
\ga\circ v(\partial_1c)= \ga\big(v(\partial_1c)\big)
\]
for any $2$--simplex $c$. Since $\ga$ is injective,
$v$ is a $1$--cocycle. 
\end{proof}
Given a $1$--cocycle $z\in \Zrm^1(K,G)$, a \emph{cross section of}
$z$  is a function $\sr: \Si_0(U)\rightarrow G$, where $U$ is an open
set of $K$, such that 
\begin{equation}
\label{Bc:7}
 z(b)\,\sr_{\partial_1b} =  \sr_{\partial_0b}, 
\qquad b\in\Si_1(U) \ .
\end{equation}
The  cross section $s$ is said to be \emph{global} whenever $U=K$. A 
reason for the terminology cross section of a $1$--cocycle 
is provided by the following 
\begin{lemma}
\label{Bc:8}
A $1$--cocycle is a $1$--coboundary 
if, and only if, it admits a global cross section.
\end{lemma}
\begin{proof}
The proof follows straightforwardly from the definition of a global cross 
section and from the definition  of 
a $1$--coboundary.
\end{proof}
\begin{oss}
\label{Bc:9}
Given a group $G$, it is very easy to define $1$--coboundaries 
of the poset $K$ with values in $G$. It is enough to assign 
an element $\sr_a\in G$ to any $0$--simplex $a$ and set 
\[
  z(b)\defi \sr_{\partial_0b}\  
              \sr^{-1}_{\partial_1b}, \qquad b\in\Si_1(K). 
\]
It is clear that $z$ is a $1$--cocycle. It is a $1$--coboundary 
because the function $\sr:\Si_0(K)\to G$ is a global cross section of $z$. 
As we shall see in the next section, the existence of 
$1$--cocycles, which are not $1$--coboundaries,
with values in a group  $G$  is equivalent to
the existence of nontrivial group homomorphisms from the 
first homotopy group of $K$ into $G$.
\end{oss}
We call \emph{the category of $1$--cocycles} with values in $G$, 
the full subcategory of $\Crm^1(K,G)$ whose set of objects 
is $\Zrm^1(K,G)$. We denote this category 
by the same symbol  $\Zrm^1(K,G)$ as used to denote the corresponding 
set of objects. It is worth observing that, given a  group homomorphism 
$\ga:G_1\to G$, by Lemma \ref{Bc:6},  
the restriction of the functor $\ga\circ$ to $\Zrm^1(K,G_1)$ defines 
a functor from $\Zrm^1(K,G_1)$ into $\Zrm^1(K,G)$.\smallskip

We interpret $1$--cocycles of  $\Zrm^1(K,G)$ as principal bundles 
over the poset $K$, having $G$ as a structure group. 
It is very easy to see which notion 
corresponds to that of an associated bundle in this framework. 
Assume that  there is an action 
$\alpha:G\times X\ni (g,x)\to \alpha(g,x)\in X$ of $G$ on a set
$X$. Consider the 
group homomorphism 
$\tilde{\alpha}: G\ni g\to \tilde{\alpha}_g\in\mathrm{Aut}(X)$  
defined as 
\[
\tilde{\alpha}_g(x)\defi \alpha(g,x),\qquad x\in X,
\]
for any $g\in G$. 
Given a $1$--cocycle $z\in\Zrm^1(K,G)$, we call the $1$-cocycle 
\begin{equation}
\label{Bc:10}
\tilde{\alpha}\circ z \in \Zrm^1(K,\mathrm{Aut}(X)),
\end{equation}
\emph{associated} with $z$,
where $\tilde{\alpha}\circ$ is the functor, associated with the 
group homomorphism $\tilde{\alpha}$, 
from the category $\Zrm^1(K,G)$ into $\Zrm^1(K,\mathrm{Aut}(X))$.
\paragraph{Homotopy and $1$--cohomology :} The 
relation between the homotopy and  the $1$--cohomology of
$K$ has been 
established in \cite{Ruz}. Here we reformulate this result 
in the language of categories. We begin by recalling some basic
properties of $1$--cocycles.  First, 
any $1$--cocycle $z\in \Zrm^1(K,G)$ 
is \emph{invariant under homotopy}. To be
precise given a pair of paths $p$ and $q$ with the same endpoints, we
have 
\begin{equation}
\label{Bc:11}
p\sim q \ \ \Rightarrow \ \ z(p) = z(q).
\end{equation}
Secondly, the following properties hold:
\begin{equation}
\label{Bc:12}
\begin{array}{lcccl}
(a) & z(\overline{p}) & = & z(p)^{-1} \ ,  &  \mbox{ for any path } p;\\
(b) & z(\sigma_0(a)) & = & e \ ,  &  \mbox{ for any $0$--simplex } a,
\end{array}
\end{equation}
Now in order 
to relate the homotopy of a poset 
to $1$--cocycles, a preliminary definition is necessary.\\ 
\indent Fix  a group $S$. Given a group $G$ 
we denote the set of group homomorphisms 
from $S$ into  $G$ by  $H(S,G)$. For any pair $\si,\si_1\in H(S,G)$   
a \emph{morphism} from $\si_1$ to $\si$ 
is an element $h$ of $G$ such that 
\begin{equation}
\label{Bc:13}
 h\ \si_1(g)  = \si(g)\  h, \qquad g\in S. 
\end{equation}
The set of morphisms from $\si_1$ to $\si$
is denoted by $(\si_1,\si)$ and there is an obvious composition rule 
between morphisms yielding a category again denoted by $H(S,G)$. 
Given a group homomorphism $\ga:G_1\to G$, there is a functor 
$\ga\circ:H(S,G_1)\to H(S,G)$ defined as 
\begin{equation}
\label{Bc:14}
\begin{array}{ll}
\ga\circ\si  \defi \ga\si    \qquad & \si \in H(S,G_1);\\
\ga\circ h \defi  \ga(h)    & h\in (\si,\si_1), \ \si,\si_1\in  H(S,G_1).   
\end{array}
\end{equation}
When $\ga$ is a group isomorphism, then $\ga\circ$ is  an isomorphism 
of categories, too. Similarly, let $S_1$  be a group
and let $\rho: S_1\rightarrow S$ be a group homomorphism. Then there is  
a functor $\circ\rho: H(S,G)\rightarrow H(S_1,G)$ defined by
\begin{equation}
\label{Bc:15}
\begin{array}{ll}
\si\circ\rho  \defi \si\rho   \qquad & \si \in H(S,G);\\
 h\circ\rho  \defi  h    & h\in (\si,\si_1), \ \si,\si_1\in  H(S,G).   
\end{array}
\end{equation}
When $\rho$ is a group isomorphism, then $\circ\rho$ is  an isomorphism 
of categories, too.\\
\indent Now, fix a base $0$--simplex $a_0$ and consider 
the category $H(\pi_1(K,a_0),G)$ associated with the first 
homotopy group of the poset. Then 
\begin{prop}
\label{Bc:16}
Given a group $G$ and  any $0$--simplex $a_0$ the categories 
$\Zrm^1(K,G)$ and $H(\pi_1(K,a_0),G)$ are equivalent.
\end{prop}
\begin{proof} 
Let us start by defining a functor 
from $\Zrm^1(K,G)$ to  $H(\pi_1(K,a_0),G)$. 
Given $z,z_1\in\Zrm^1(K,G)$ and $\fr\in(z_1,z)$,
define
\[
   \begin{array}{ll}
    F(z)([p]) \defi z(p), &  [p]\in\pi_1(K,a_0);\\
    F(\fr) \defi \fr_{a_0}.  
   \end{array}
\]
$F(z)$ is well defined since $1$--cocycles are
homotopy invariant. Moreover, it is easy to see 
by (\ref{Bc:12}) that $F(z)$ is a group homomorphism from 
$\pi_1(K,a_0)$ into $G$. Note that 
\[
 \fr_{a_0}\  F(z_1)([p]) =  \fr_{a_0} \  z_1(p) = 
 z(p)\  \fr_{a_0} = F(z)([p])\  \fr_{a_0},
\]
hence $F(\fr)\in (F(z_1), F(z))$. So $F$ is well defined and 
easily shown to be a covariant functor. 
To define a functor $C$ in the other direction, let us choose 
a path $p_a$  from $a_0$ to $a$, for any $a\in\Si_0(K)$.
In particular we set $p_{a_0} = \sigma_0(a_0)$. 
Given $\si\in H(\pi_1(K,a_0),G)$  and  
$h\in(\si_1,\si)$, define 
\[
   \begin{array}{ll}
    C(\si)(b) \defi
    \si([\overline{p}_{\partial_0b}*b*p_{\partial_1b}])\ , & 
       b\in\Si_1(K); \\
    C(h) \defi \mathrm{c}(h),
   \end{array}
\]
where $\mathrm{c}(h):\Si_0(K)\to G$ is the constant function 
taking the value  $h$ for any $a\in\Si_0(K)$. 
It can be easily shown that $C$ is a covariant functor. 
Concerning the equivalence, note that 
\[
 (F\cdot C)(\si)([p]) = C(\si)(p) = \si([\overline{\sigma_0(a_0)}*p* 
\sigma_0(a_0)])    = \si([p]), 
\]
and that 
\[
 (F\cdot C)(h) = F(\mathrm{c}(h))  =  h.
\]
Hence $F\cdot C = \mathrm{id}_{H(\pi_1(K,a_0),G)}$.  
Conversely, given a $1$--simplex $b$ we have 
\[
(C\cdot F)(z)(b)= 
F(z)([\overline{p_{\partial_0b}}*b*p_{\partial_1b}]) = 
z(p_{\partial_0b})^{-1}\  z(b)\   z(p_{\partial_1b}) ,
\]
and  given a $0$--simplex $a$ we have 
\[
 (C\cdot F)(\fr) =  C (\fr_{a_0}) = \mathrm{c}(\fr_{a_0})
\]
Define $u(z)_a\defi z(p_a)$ for $a\in\Si_0(K)$. It can be easily seen 
that the mapping $\Zrm^1(K,G)\ni z\to u(z)$ defines 
a natural isomorphism between the functor $C\cdot F$ and 
the functor $\mathrm{id}_{\Zrm^1(K,G)}$. 
\end{proof}
Observe in particular that the group homomorphism corresponding 
to the trivial $1$--cocycle $\imath$ is the trivial 
one, namely $\si([p]) = e$ for any $[p]\in \pi_1(K,a_0)$.  
Hence, a $1$--cocycle of $\Zrm^1(K,G)$ is a $1$--coboundary if, and only 
if, the corresponding group homomorphism $F(z)$ is equivalent to the 
trivial one. In particular if $K$ is simply 
connected,
then $\Zrm^1(K,G)=\Brm^1(K,G)$ for  any group $G$.\\ 
\indent The existence of $1$--cocycles, 
which are not $1$--coboundaries, relies, in particular, 
on the following corollary 
\begin{cor}
\label{Bc:17}
Let $M$ be a nonempty, Hausdorff and  arcwise connected 
topological space which admits a base
for the topology consisting of  arcwise and simply connected
subsets of $M$. Let $K$ denote the poset formed by such a  base 
ordered under inclusion $\subseteq$. Then 
\[
H(\pi_1(M,x_0),G)\cong  H(\pi_1(K,a_0),G)\cong \Zrm^1(K,G) \ , 
\]
for any $x_0\in M$ and $a_0\in \Si_0(K)$ with $x_0\in a_0$, where the symbol 
$\cong$ means equivalence of categories. 
\end{cor}
\begin{proof}
$\pi_1(M,x_0)$ is isomorphic to $\pi_1(K,a_0)$ (see in Section \ref{A}). 
As observed at the beginning of this section, this entails 
that  $H(\pi_1(M,x_0),G)$  and 
$H(\pi_1(K,a_0),G)$ are isomorphic categories. 
Therefore the proof follows by 
Proposition \ref{Bc:16}.
\end{proof}
Let $M$ be  a nonsimply connected topological 
space and let $K$ be a basis for the topology of $M$ 
as defined in the statement  of Corollary \ref{Bc:17}.
Then to any nontrivial group homomorphism in   $H(\pi_1(M,x_0),G)$
there corresponds a $1$--cocycle of $\Zrm^1(K,G)$ which is not a 
$1$--coboundary.
\section{Connections}
\label{C} 
This section is entirely devoted to studying connections 
and related notions like the curvature, holonomy group 
and the central connections.  
We will show how connections and $1$--cocycles are related, thus
allowing one to interpret a $1$--cocycle as a principal bundle 
and a connection $1$-cochain  as the connection of this principal bundle. 
We will prove the existence of nonflat connections,
a ``poset'' version of the Ambrose-Singer Theorem, and that 
to any flat connection with values in $G$, there corresponds 
a group homomorphism from the fundamental group of the poset into $G$.
\subsection{Connections and curvature}
\label{Ca} 
We now give the definition of a connection of a poset 
with values in a group. To this end, recall the definition of the set 
$\Si^{\inf}_n(K)$ of 
inflating $n$--simplices (see Section \ref{A}).
\begin{df}
\label{Ca:1}
A $1$--cochain  $u$ of $\Crm^1(K,G)$ is said to be a 
\textbf{connection} if it satisfies the following properties:
\begin{itemize}
\item[(i)] $u(\overline{b})= u(b)^{-1}$ for any $b\in\Sigma_1(K)$; 
\item[(ii)] $u(\partial_0c)\,  u(\partial_2c) = u(\partial_1c)$,
for any $c\in \Si^{\inf}_2(K)$. 
\end{itemize}
We denote  the set of connection $1$--cochains with values 
in  $G$ by $\Urm^1(K,G)$.
\end{df}
This definition of a connection  is  related to the notion of the link 
operator in a lattice gauge theory  
(\cite{Cre}) and to the notion of a  generalized connection in 
loop quantum gravity (\cite{Bae,Lew}).  
Both the link operator and the generalized connection can be seen  
as a mapping $A$ which associates  an element 
$A(e)$ of a group $G$ to any oriented edge $e$ of a graph $\alpha$, 
and enjoying the following properties 
\begin{equation}
\label{Ca:1a}
 A(\overline{e})= A(e)^{-1}, \ \ \ A(e_2* e_1) = A(e_2)\,  A(e_1), \qquad 
\end{equation}
where, $\overline{e}$ is the reverse of the edge $e$; 
$e_2* e_1$ is the composition of the edges $e_1$, $e_2$
obtained by composing the end of $e_1$ with the beginning of $e_2$. 
Now, observe  that to any poset $K$  there corresponds 
an oriented  graph $\alpha(K)$ whose  set of vertices  is 
$\Si_0(K)$, and whose  set of edges is $\Si_1(K)$. 
Then, by property (i) of the above definition and 
property (\ref{Bc:1}), any  connection $u\in\Urm^1(K,G)$ 
defines a mapping from the edges of  $\alpha(K)$ to  $G$ satisfying 
(\ref{Ca:1a}). The new feature of our definition of connection,  
is to  require  property (ii) in  Definition \ref{Ca:1}, thus involving
the poset structure. The motivation for this property 
will become clear in the next section: thanks to this property  
any connection $u$ can be seen as a connection on the principal 
bundle described by  a 1-cocycle (see Theorem \ref{Cb:1}). \smallskip
 
Let us now observe that any $1$--cocycle is a connection. Furthermore,
if $u$ is a connection  then so is 
any equivalent $1$--cochain (the proof is similar to the proof 
of the same property for $1$--cocycles, see Section \ref{Bb}). 
\begin{lemma}
\label{Ca:2}
Let $\gamma:G_1\to G$ be a group homomorphism. Given 
$v\in \Crm^1(K,G_1)$ consider $\ga\circ v\in\Crm^1(K,G)$. 
Then: if $v$ is a connection 
then $\ga\circ v$ is a connection; the converse holds  
if $\ga$ is injective. 
\end{lemma}
\begin{proof}
Clearly, if $v$ is a connection so is $\ga\circ v$. Conversely, 
assume that $\ga$ is injective and that $\ga\circ v$ is a connection. 
If $c\in\Si^{\inf}_2(K)$, then 
\[
\ga\big(v(\partial_0c)\,  v(\partial_2c)\big)  = 
\ga\circ v(\partial_0c)\,  \ga\circ v(\partial_2c)
  =  \ga\circ v(\partial_1c) = \ga\big(v(\partial_1c)\big),
\]
hence $v(\partial_0c)\, 
v(\partial_2c)= v(\partial_1c)$,
since $\ga$ is injective. Furthermore, for any $1$--simplex $b$
we have 
\[
\ga(v(\overline{b})) =  \ga\circ v(\overline{b})= 
                          (\ga\circ v(b))^{-1} = \ga(v(b)^{-1}).
\]
So, as $\ga$ is injective, we have  $v(\overline{b})=v(b)^{-1}$,
and this entails that $v$ is a connection. 
\end{proof}
\begin{lemma}
\label{Ca:3}
Given $u\in\Urm^1(K,G)$, then \\
(a) $ u(\sigma_0(a_0)) = e$ for any $a_0\in\Si_0(K)$,\\
(b) $ u(b) =  u(b_1)$ for any $b,b_1\in\Si^{\inf}_1(K)$ with 
    $\partial b = \partial b_1$.
\end{lemma}
\begin{proof}
(a) Since a degenerate $1$--simplex is an inflating $1$--simplex,  
by Definition \ref{Ca:1}(ii) we have 
\[
  u(\sigma_0(a_0))\   u(\sigma_0(a_0)) =  u(\sigma_0(a_0)) 
 \ \iff \  u(\sigma_0(a_0)) = e. 
\]
(b) Given $b_1$ and $b$ as in the statement. Since $b,b_1$
are inflating $1$--simplices, 
the $1$--simplex $b_0$ defined as 
$\partial b_0 \defi\partial b$, and  $|b_0|  \defi \partial_0b$,  
is inflating. Moreover $\partial b_0=\partial b=\partial b_1$
and $|b_0|\subseteq |b|, |b_1|$. As 
$|b_0|\subseteq |b|$, the $2$--simplex $c$
defined by 
\[
 \partial_0c \defi \sigma_0(\partial_0b), \ \ 
 \partial_1c \defi  b, \ \ 
 \partial_2c \defi  b_0, \ \  |c| \defi |b| 
\]
is an inflating $2$--simplex. By Definition \ref{Ca:1}(ii) 
and by (a) we have 
\[
 u(b_0) =  u(\si_0(\partial_0b))\,   u(b_0) = 
 u(\partial_0c)\,    u(\partial_2c) = 
 u(\partial_1c) =  u(b) \ .
\]
The same reasoning leads to $u(b_0)= u(b_1)$, hence 
$u(b)=  u(b_1)$.
\end{proof}
In words, this lemma says that connections act 
trivially on degenerate $1$--simplices, and that 
their values do not depend on the support of the 
inflating $1$--simplices.\smallskip  

We call the full subcategory of $\Crm^1(K,G)$ whose set of objects 
is $\Urm^1(K,G)$ \emph{the category of connection $1$--cochains} 
with values in $G$. It will be  denoted 
by the same symbol  $\Urm^1(K,G)$ as used to denote the corresponding 
set of objects. Note that  $\Zrm^1(K,G)$ is a full subcategory 
of $\Urm^1(K,G)$. Furthermore, if  $\ga:G_1\to G$ is a  group homomorphism,  
by Lemma \ref{Ca:2},  the restriction of the functor 
$\ga\circ$ to $\Urm^1(K,G_1)$ defines 
a functor from $\Urm^1(K,G_1)$ into $\Urm^1(K,G)$.\smallskip

As observed,  any  $1$--cocycle is a connection. 
The converse does not hold, in general, and the 
obstruction is  a $2$--coboundary.
\begin{df}
\label{Ca:4}
The \textbf{curvature } of  a connection $u\in\Urm^1(K,G)$ is the 
$2$--coboundary  $W_u\defi \dr u\in\Brm^2(K,G)$. 
Explicitly, by using relation (\ref{Bb:11})  we have 
\[
\begin{array}{rcll}
 (W_u)_1(b) &=& \big( e,\mathrm{ad}(u(b))\big), & b\in\Si_1(K), \\
 (W_u)_2(c)&=& \big(w_u(c), \mathrm{ad}(u(\partial_1c)\big) & c\in\Si_2(K),
\end{array}
\] 
where 
$w_u:\Si_2(K)\to G$ defined as 
\[
w_u(c) \defi u(\partial_0c)\ u(\partial_2c)\ u(\partial_1c)^{-1}, \qquad c\in\Si_2(K). 
\]
A connection $ u\in\Urm^1(K,G)$ is said to be \textbf{flat} whenever 
its curvature is trivial i.e. $W_u\in(2G)_1$ or, equivalently, 
if $w_u(c)=e$ for any $2$--simplex  $c$.
\end{df}
We now draw some consequences of our definition of the curvature
of a connection  and  point out the relations of this notion 
to  the corresponding one  in the
theory of principal bundles. \smallskip 
 
\emph{First},  note that a connection $u$
is flat if, and only if, $u$ is a $1$--cocycle. Then,  
as an immediate consequence of Proposition \ref{Bc:16}, 
we have  a poset version of a classical result of the theory 
of principal bundles \cite{KN, DK}.
\begin{cor}
\label{Ca:5}
There is, up to equivalence, a 1-1 correspondence 
between flat connections of $K$ with values in $G$ and 
group homomorphisms from $\pi_1(K)$ into $G$.   
\end{cor}
The existence of nonflat connections 
will be shown in Section \ref{Cc} where examples will be given.\smallskip 

\emph{Secondly}, in a principal bundle  the curvature form  
is the covariant exterior derivative of a connection form, namely 
the $2$--form with values in the Lie algebra of the group, 
obtained by taking the exterior derivative 
of the connection form and evaluating this on the horizontal
components of pairs of vectors of the tangent space (see \cite{KN}). 
Although, no differential structure is present in our approach,  
but $W_u$ encodes this type of information. 
In fact, given a connection $u$, if we interpret 
$u(p)$ as the horizontal lift of a path $p$, then the equation
\begin{equation}
\label{Ca:6}
 w_u(c)\, u(\partial_1c)  = 
          u(\partial_0c* \partial_2c)\, w_u(c), \qquad 
   c\in\Si_2(K),  
\end{equation}
may be understood as saying  that $w_u(c)$ intertwines 
the horizontal lift of the path $\partial_1c$  and that 
of the path $\partial_0c* \partial_2c$. \smallskip

\emph{Thirdly}, the structure equation of the curvature form (see \cite{KN}) 
says that the curvature equals the exterior derivative  
of the connection form plus the commutator of the connection form. 
Notice that the second component $(W_u)_2$ of the curvature 
can be rewritten as 
\begin{equation}
\label{Ca:7}
 (W_u)_2(c) = \big(w_u(c),\mathrm{ad}(w_u(c))^{-1}\big)\times 
              \big(e,\mathrm{ad}(u(\partial_0c)u(\partial_2c))\big), 
\end{equation}
for any $2$--simplex $c$,  
where $\times$ is the composition (\ref{Bb:2}) of the $2$--category 
$2G$. This, equation represents, in our formalism, 
the structure equation of the curvature with 
$\big(w_u(c),\mathrm{ad}(w_u(c)^{-1})\big)$ in place of the exterior 
derivative, and 
$\big(\io,\mathrm{ad}(u(\partial_0c)u(\partial_2c))\big)$  in 
place of the commutator of the connection form. \smallskip

\emph{Fourthly}, as a consequence of  Lemma \ref{Bb:13} we have that 
$W_u$ is a $2$--cocycle. The $2$--cocycle identity 
that is $\dr W_u\in (3G)_2$ or, equivalently,
\begin{equation}
 w_u(\partial_0d) \ w_u(\partial_2d)  = 
 \mathrm{ad}(u(\partial_{01}d))\big(w_u(\partial_3d)\big) \ 
                           w_u(\partial_1d), \qquad \Si_3(K),
\end{equation}
corresponds, in our framework, to  the Bianchi identity.\smallskip

We conclude with the following result.
\begin{lemma}
\label{Ca:8}
For any connection $u$ the following assertions hold: \\
(a) $w_u(c^{\sst{(01)}}) =  w_u(c)^{-1}$ for any $2$--simplex $c$; \\
(b) $w_u(c)=e$ if $c$ is either a degenerate or an inflating
$2$-simplex. 
\end{lemma}
\begin{proof}
(a) follows directly from the definition of $c^{\sst{(01)}}$, see
Section \ref{A}.
(b) If $c$ is an inflating $2$--simplex, then
$w_u(c)=e$ because of Definition \ref{Ca:1}(ii).
Given a $1$--simplex $b$, then 
\begin{align*}
w_u(\si_0(b)) & = u(\partial_0\si_0(b)) \  
                 u(\partial_2\si_0(b)) \  u(\partial_1\si_0(b))^{-1} \\
              & =  u(b) \  
                 u(\si_0(\partial_1b)) \   u(b)^{-1} = e, 
\end{align*}
because by Lemma \ref{Ca:3}(a) we have that
$u(\si_0(\partial_1b))=e$.
Analogously we have that $w_u(\si_1(b))=e$.
\end{proof}
In words,  statement (b) asserts that 
the curvature of a connection  is trivial when restricted 
to inflating simplices.
\begin{oss}
\label{Ca:9}
It is worth pointing out  some analogies 
between the theory of connections, as presented in this paper, 
and that developed in synthetic geometry 
by A. Kock \cite{Kock1}, and 
in algebraic topology  by   L. Breen and W. Messing \cite{BM1}. 
The contact point with our approach resides 
in the fact that both of the other approaches make use of a 
combinatorial notion of differential forms taking values in a group
$G$. So in both  cases connections turn out to be 
combinatorial $1$--forms. Concerning the curvature, 
the  definition of $W_u$ is formally the same as
the definition of curving data  given 
in \cite{BM1}, since this is the $2$--coboundary of a connection, 
taking values in a $2$--category associated with $G$. Whereas, 
in \cite{Kock1} the curvature is the $2$--coboundary 
of a connection,  taking values 
in $G$, and  is formally the same as $w_u$.  
The only difference to  these other two approaches 
is that in our case $w_u$ is not invariant under oriented equivalence of 
$2$-simplices (examples 
of connections having this feature will be given at the end 
of  Section \ref{Cc}).  
\end{oss}
%
\subsection{The cocycle induced by a connection}
\label{Cb}
We analyze  the relation between 
connections and $1$--cocycles more deeply. The main result is that 
to any connection there corresponds a unique $1$--cocycle. 
This, on the one hand, confirms the interpretation of $1$--cocycles
as principal bundles. On the other hand this result will allow us to 
construct
examples of nonflat connections in the next section.
\begin{teo}
\label{Cb:1}
For any $u\in \Urm^1(K,G)$, there  exists 
a unique $1$--cocycle $z\in\Zrm^1(K,G)$ such that 
$u(b) = z(b)$ for any $1$--simplex  $b\in\Si^{\inf}_1(K)$.
\end{teo} 
\begin{proof}
Within this proof we adopt 
the following notation: for any $\dc\in U_a\cap U_{a_1}$ the 3--tuple 
$(\dc;a,a_1)$ denotes the $1$--simplex  with support $\dc$, 
$0$--boundary $a$ and $1$--boundary $a_1$.  Consider the open 
set $U_a$  (\ref{A:10}) of the fundamental covering of $K$, and define 
\begin{equation}
\label{Cb:2}
 z_{a_1,a}(\dc) \defi u(\dc;\dc,a_1)^{-1} \, u(\dc;\dc, a), 
    \qquad \dc\in U_a\cap U_{a_1}.
\end{equation}
So, we have a family of functions 
$z_{a_1,a}: U_a\cap U_{a_1}\to G$. 
Let  $\dc_1\subseteq \dc$, with $\dc,\dc_1\in U_a\cap U_{a_1}$.
By using the defining properties of connection 
we have 
\begin{align*}
z_{a_1,a}(\dc) & = u(\dc;\dc,a_1)^{-1} \  u(\dc;\dc, a) \\
  & =  \big(u(\dc;\dc,\dc_1)\   u(\dc_1;\dc_1,a_1)\big)^{-1}
    \  u(\dc;\dc,\dc_1) \  u(\dc_1;\dc_1, a) \\
   & = u(\dc_1;\dc_1,a_1)^{-1} \  u(\dc_1;\dc_1, a) \\
  & = z_{a_1,a}(\dc_1).
\end{align*}
If $a\subseteq a_1$ and $\dc\in U_{a_1}$,  then 
\begin{align*}
 z_{a_1,a}(\dc) & = u(\dc;\dc,a_1)^{-1} \,  u(\dc;\dc, a) 
  = u(\dc;\dc,a_1)^{-1}\,  u(\dc;\dc, a_1) \,  u(\dc;a_1, a) \\
 & =u(\dc;a_1, a)
\end{align*}
Assume that $\dc\in U_a\cap U_{a_1}\cap U_{a_2}$. Then 
\begin{align*}
z_{a_2,a_1}(\dc)\,  z_{a_1,a}(\dc) 
 & = u(\dc;\dc,a_2)^{-1} \,  u(\dc;\dc, a_1)\,  
   u(\dc;\dc,a_1)^{-1} \,  u(\dc;\dc, a)\\
 & = u(\dc;\dc,a_2)^{-1} \,  u(\dc;\dc, a)
  = z_{a_2,a}(\dc).
\end{align*}
Summing up, to a  connection $ u$  corresponds a family 
of functions $z_{a_1,a}: U_a\cap U_{a_1}\to G$  satisfying the following 
properties: 
\begin{equation}
\label{Cb:3}
\begin{array}{rll}
(i) & z_{a_1,a}(\dc) = z_{a_1,a}(\dc_1), \ \ \ & 
 \dc_1\subseteq \dc,  \ \ \dc,\dc_1\in U_a\cap U_{a_1}; \\
(ii) &  z_{a_1,a}(\dc) =u(\dc;a_1, a), & 
           \dc\in U_{a_1}, \ \ a\subseteq a_1; \\
(iii) & z_{a_2,a_1}(\dc)\,  z_{a_1,a}(\dc) 
 = z_{a_2,a}(\dc),  & \dc\in U_a\cap U_{a_1}\cap U_{a_2}.
\end{array}
\end{equation}
Now, note that since $\partial_0b,\partial_1b\leq |b|$, we have 
that $|b|\in U_{\partial_0b}\cap U_{\partial_1b}$. Hence we can define
\begin{equation}
\label{Cb:4}
z(b)\defi z_{\partial_0b,\partial_1b}(|b|), \qquad b\in\Si_1(K).  
\end{equation}
Given a $2$--simplex $c$. By using  properties (i)--(iii) 
we have 
\begin{align*}
 z(\partial_0c)\  z(\partial_2c) 
& =  z_{\partial_{00}c,\partial_{10}c}(|\partial_0c|)\,  
     z_{\partial_{02}c,\partial_{12}c}(|\partial_2c|) \\
& =  z_{\partial_{00}c,\partial_{10}c}(|c|)\,  
     z_{\partial_{02}c,\partial_{12}c}(|c|) \\
& =  z_{\partial_{01}c,\partial_{10}c}(|c|)\,  
     z_{\partial_{10}c,\partial_{11}c}(|c|) \\
& =  z_{\partial_{01}c,\partial_{11}c}(|c|) \\
& =  z(\partial_{1}c). 
\end{align*}
Hence $z$ is $1$-cocycle. 
Moreover, if $b$ is an inflating $1$--simplex, then  
\begin{align*}
z(b) & = z_{\partial_0b,\partial_1b}(|b|) =
       u(|b|;|b|,\partial_0b)^{-1} \,  u(|b|;|b|, \partial_1b) \\
      & = u(|b|;|b|,\partial_0b)^{-1} \,  
          u(|b|;|b|, \partial_0b)\,  u(|b|;\partial_0b, \partial_1b)
       = u(|b|;\partial_0b, \partial_1b) \\
       & = u(b).
\end{align*}
$z$ is clearly the unique $1$--cocycle with $z(b)=u(b)$ for 
$b\in\Sigma_1^{\inf}(K)$.
\end{proof}
On the basis of this theorem, we can introduce the following
definition. 
\begin{df}
\label{Cb:5}
A  connection 
$u\in\Urm^1(K,G)$ is said to \textbf{induce} the 
$1$--cocycle $z\in\Zrm^1(K,G)$ whenever 
\[
 u(b) = z(b), \qquad  b\in\Si^{\inf}_1(K) \ .
\]
We denote  the set of connections of 
$\Urm^1(K,G)$ inducing 
the $1$--cocycle $z$ by  $\Urm^1(K,z)$.
\end{df}
The geometrical meaning of $\Urm^1(K,z)$ is the
following: just as a $1$--cocycle $z$ stands for a principal bundle over $K$ so
the set of connections $\Urm^1(K,z)$ stands for the set of 
connections on that principal bundle.
Theorem \ref{Cb:1} says that the set of connections 
with values in $G$ is partitioned as 
\begin{equation}
\label{Cb:6}
 \Urm^1(K,G) = \dot{\cup}\big\{ \Urm^1(K,z) \ | \
 z\in\Zrm^1(K,G)\big\} \ 
\end{equation}
where the symbol $\dot{\cup}$ means disjoint union. 
\begin{lemma}
\label{Cb:7}
Given  $ z_1,z\in\Zrm^1(K,G)$,  let 
$u_1\in U^1(K, z_1)$ and $u\in \Urm^1(K,z)$. Then 
$(u_1,u)\subseteq(z_1,z)$. In particular 
if $u_1\cong u$, then $z_1\cong z$.
\end{lemma}
\begin{proof}
(a) By equations (\ref{Cb:2}) and (\ref{Cb:4}), we have 
\[
 z(\dc;a_1,a) = u(\dc;\dc,a_1)^{-1}\  u(\dc; \dc,a), 
\]
for any $1$--simplex $(\dc;a_1,a)$.
The same holds for $z_1$ and $u_1$.
Given $\fr\in(u_1,u)$, we have 
\begin{align*}
\fr_{a_1}\,  z_1(\dc;a_1,a) & =
\fr_{a_1}\,  u_1(\dc;a_1,\dc)\,  u_1(\dc;\dc,a) \\
&  = u(\dc;a_1,\dc)\,  \fr_{\dc} \,  
     u_1(\dc;\dc,a) \\
& = u(\dc;\dc,a_1)^{-1} \,  u(\dc; \dc,a) \,  \fr_{a} 
 = z(\dc;a_1,a) \,  \fr_{a},
 \end{align*} 
where we have use the fact that  $(\dc;a_1,\dc)$ is the reverse
of $(\dc;\dc,a_1)$. Hence $\fr\in (z_1,z)$, 
and this completes the proof. 
\end{proof}
Now, given a $1$-cocycle $z\in\Zrm^1(K,G)$, we call the category of
\emph{connections inducing} $z$, the full subcategory 
of $\Urm^1(K,G)$ whose objects belong to $\Urm^1(K,z)$. 
As it is customary in this paper, we denote this category 
by the same symbol $\Urm^1(K,z)$ as used to denote the corresponding
set of objects. 
\begin{lemma}
\label{Cb:8}
Let $z\in\Zrm^1(K,G_1)$ and let $\ga:G_1\to
G$ be an injective group homomorphism. Then, the functor 
$\ga\circ : \Urm^1(K,z)\to \Urm^1(K,\ga\circ z)$
is injective and faithful.
\end{lemma}
\begin{proof}
Given $u\in \Urm^1(K z)$, it is easy to see that 
$\ga\circ u\in \Urm^1(K,\ga\circ z)$.  Clearly, 
as $\ga$ is injective, the functor  $\ga\circ $ is injective and faithful. 
\end{proof}
We note the following simple result.
\begin{lemma}
\label{Cb:9}
If $z_1\cong z$,  then the categories $\Urm^1(K,z_1)$ and
$\Urm^1(K,z)$  are equivalent.
\end{lemma}
Assume that  $K$ is simply connected. In this case  any $1$-cocycle 
is a $1$--coboundary (see Section \ref{Bc}). Then the category 
$\Urm^1(K,z)$ is equivalent to  $\Urm^1(K,\imath)$  for any  $z\in\Zrm^1(K,G)$.
\subsection{Central connections}
\label{C:x}
We now briefly study the family  of central connections, whose 
main feature, as we will show below, is that any such  connection can be 
uniquely decomposed as the product of the induced cocycle 
by a suitable connection taking values in the centre of the group. 
\begin{df}
\label{Cx:1}
A  connection $u\in\Urm^1(K,G)$ is said to be 
a \textbf{central connection} whenever 
the component $w_u$ of the curvature $W_u$ takes values in the centre 
$\mathcal{Z}(G)$.
We denote the set of central  connections by 
$\Urm^1_{\mathcal{Z}}(K,G)$. 
\end{df}
Let us start to analyze the properties of central connections. 
Clearly $1$--cocycles are central
connections. However, the  main property that can be directly 
deduced from the above definition is that the component 
$w_u$ of the curvature $W_u$ of
a central connection $u$ is invariant under oriented 
equivalence of $2$--simplices. 
In fact by the definition of $w_u$,   
it is easily seen that 
\begin{equation}
\label{Cx:2}
 w_u(c) = w_u(c_1) = w_u(c_2)^{-1}, \qquad c_1\in[c],
c_2\in\overline{[c]},
\end{equation}
for any $2$--simplex $c$, 
where $[c]$ and $\overline{[c]}$ are, respectively,  the classes of 
$2$--simplices having the same and the reversed 
orientation of $c$. 
\begin{prop}
\label{Cx:3}
A connection   $u$  of $\Urm^1(K,G)$ is central 
if, and only if, it can be uniquely decomposed as
\[
   u(b) = z_u(b) \, \chi_u(b), \qquad b\in\Si_1(K), 
\] 
where $z_u\in \Zrm^1(K,G)$, and 
$\chi_u\in \Urm^1(K,\imath)$ with values in $\mathcal{Z}(G)$. 
\end{prop}
\begin{proof}
$(\Leftarrow)$ Assume that  a connection $u$ admits a decomposition as
in the statement. Since $\chi_u$ takes values in the centre, 
so does $w_u$. Furthermore, 
since $\chi_u\in\Urm^1(K,\imath)$ then $\chi_u(b)= e$ for any 
inflating $1$--simplex $b$. This entails that $z_u$ is nothing but
the $1$--cocycle induced by $u$. This is enough for uniqueness. 
$(\Rightarrow)$ Assume that $u$ is central. For any  $1$--simplex $b$ let 
$c_b$ denote the $2$--simplex defined as 
\[
  \partial_1c_b \defi  b, \ \   \partial_0c_b \defi  
  (|b|,\partial_1b, |b|), \ \ 
  \partial_2c_b\defi  (|b|,|b|, \partial_0b), \ \ |c_b|\defi |b|.
\] 
As  $w_u$ takes values in $\mathcal{Z}(G)$ we have
\[
 u(b)\,  w_u(c_b) = u(\partial_1c_b)\, w_u(c_b)  = 
  u(\partial_0c_b)\, u(\partial_2c_b) = z_u(b),  
\]
where the latter identity is a consequence of 
the fact that $u(\partial_0c_b)\, u(\partial_2c_b)$ is nothing but
the definition (\ref{Cb:2}) of the 1-cocycle induced by $u$. 
Now, define 
\[
\chi_u(b)\defi w_u(c_b), \qquad b\in\Si_1(K).
\]
Since $\chi_u(b) = u(b)\, z_u(b)^{-1}$, one can easily deduce that 
$\chi_u\in \Urm^1(K,\imath)$, and this completes the proof.
\end{proof}
As a consequence of this result 
the set $\Urm^1_{\mathcal{Z}}(K,z)$ of  central connections 
inducing the $1$--cocycle $z$, has a the structure of an Abelian group. 
In fact,  given 
$u,u_1\in \Urm^1_{\mathcal{Z}}(K,z)$, define 
\begin{equation}
\label{Cx:4}
 u\star_z u_1 (b) \defi u(b) \ z(b)^{-1}\ u_1(b), \qquad b\in\Si_1(K). 
\end{equation}
By Proposition \ref{Cx:3}, we have 
$u\star_z u_1 (b) = z(b)\ \chi_u(b)\ \chi_{u_1}(b)$ for any 
$1$--simplex $b$. This entails that 
\[
u\star_z u_1=u_1\star_z u \ \mbox{ and } \ 
u\star_z u_1\in \Urm^1_{\mathcal{Z}}(K,z). 
\]
Finally, it can be easily 
seen that $\Urm^1_{\mathcal{Z}}(K,z)$ with $\star_z$ is an Abelian  
group whose identity is $z$, and such that the inverse 
of a connection $u$ is the connection defined as 
$z(b)\ \chi_u(b)^{-1}$ for any $1$--simplex $b$.\smallskip 

Finally, in Section \ref{Ca} we pointed out the analogy between 
equation (\ref{Ca:7}) and the structure equation  of  
the curvature of a connection in a principal bundle. 
This analogy is stronger for a central connection $u$
since we have 
\begin{equation}
\label{Cx:5}
\begin{array}{lcl}
  (W_u)_2(c)
  & = & \big(w_u(c),\io\big)\times 
              \big(e,\mathrm{ad}(u(\partial_0c)u(\partial_2c))\big) \\[3pt]
           &=&
              \big(e,\mathrm{ad}(u(\partial_0c)u(\partial_2c))\big) 
              \times \big(w_u(c),\io \big),
\end{array}
\end{equation}
for any $2$--simplex $c$. Hence, as for  principal bundles, 
equation (\ref{Ca:7})  for a central connection is symmetric 
with respect to the interchange of the two factors.  
\subsection{Existence of nonflat connections}
\label{Cc}
We investigate the existence of nonflat connections. 
As a first step, we show  that there is  a very particular class of 
posets which not admitting  
nonflat connections. \smallskip

Recall that a poset $K$ is said to be \emph{totally ordered} 
whenever for any pair  $\dc,\dc_1\in K$ either 
$\dc\leq \dc_1$ or $\dc_1\leq\dc$.  Clearly, a totally ordered poset 
is directed and, consequently, pathwise connected (it is also simply
connected, see 
Section \ref{A}). 
\begin{cor}
\label{Cc:1}
If $K$ is totally ordered, any connection  is flat. 
\end{cor}
\begin{proof}
If $K$ is totally ordered  and $b$ is any $1$--simplex  either $b$ or 
$\overline{b}$ is an  inflating $1$--simplex. 
Hence, by Theorem \ref{Cb:1}  any connection coincides with the 
associated $1$--cocycle. 
\end{proof}
Another obvious situation where  nonflat connections do not 
exist is when the group  of  coefficients $G$ is trivial, i.e. $G=e$. 
Two observations on these results are in order. First,  Corollary 
\ref{Cc:1} cannot be directly deduced 
from the definition of a connection.  Secondly, 
as explained earlier, these two situations never arise in the 
applications we have in mind.\\
\indent Now, our purpose is to show that except when the poset 
is totally directed or  the group of coefficients is trivial, 
nonflat connections always exist. Let us starting by the following 
\begin{lemma}
\label{Cc:2}
Assume that there exists a $1$--cochain 
$v\in\Crm^1(K,G)$ such that $v(b) = e=v(\overline{b})$ for any 
inflating $1$--simplex $b$. Then,  for any $1$--cocycle 
$z\in\Zrm^1(K,G)$ the $1$--cochain  $v(z)$ defined as 
\begin{equation}
\label{Cc:3}
 v(z)(b) \defi v(\overline{b})^{-1} \,  z(b)\,  v(b),
            \qquad b\in\Si_1(K), 
\end{equation}
is a connection  inducing $z$.
\end{lemma}
\begin{proof}
By the definition of  $v$ for any inflating $1$--simplex $b$ we have that 
\[
  v(z)(b) = v(\overline{b})^{-1} \,  z(b)\,  v(b) = 
         e \,  z(b)\,   e  = z(b) \ .
\]
This, in particular, entails that $v(z)$ satisfies property 
$(ii)$ of the definition of connections.
For any $1$--simplex $b$ we have 
\begin{align*}
   v(z)(\overline{b}) & =  v(\overline{\overline{b}})^{-1} \,   
                       z(\overline{b})\,  v(\overline{b})  
       = v(b)^{-1} \,   z(b)^{-1}\,  v(\overline{b})  \\
      & =   \left(v(\overline{b})^{-1}\,  z(b)\, 
           v(b)\right)^{-1} =  v(z)(b)^{-1}.
\end{align*}
Hence $v(z)\in\Urm^1(K,z)$.
\end{proof}
It is very easy to prove the existence of elements of $\Crm^1(K,G)$ 
satisfying the properties of the statement. For instance, 
given a  $1$--simplex $b$, define 
\begin{equation}
\label{Cc:4} 
v(b)\defi \left\{
\begin{array}{ll}
 e & b \mbox{ or } \overline{b} \in\Si^{\inf}_1(K) \\ 
 g(b) &  \mbox{ otherwise },
\end{array}\right.
\end{equation}
where $g(b)$ is some element of the group $G$. So
$v$ is a $1$--cochain satisfying the relation 
$v(b) = e=v(\overline{b})$ for any inflating $1$--simplex $b$. \\ 
\indent Now, assume that  $K$ is a pathwise connected but not
totally directed poset. Let $G$ be a nontrival group.  
Let $v\in\Crm^1(K,G)$ be defined by (\ref{Cc:4}), and let 
$z\in\Zrm^1(K,G)$. Consider the connection $v(z)\in\Urm^1(K,z)$.
We want to find conditions on $v$ implying that $v(z)$ is not flat.\\ 
\indent  As $v(z)\in \Urm^1(K,z)$,   Theorem \ref{Cb:1}  
says that if $v(z)$ is flat then $v(z)=z$. Hence, 
$v(z)$ is not flat if, and only if, it differs from 
$z$ on a $1$--simplex $b$ such that both $b$ and $\overline{b}$ are 
noninflating. Then 
\begin{align*}
v(z)(b) \ne z(b) & \iff v(\overline{b})^{-1} \,  z(b)\,  v(b) \ne z(b)\\
              & \iff z(b)\,  v(b) \ne v(\overline{b})\,  z(b)\\
              & \iff z(b)\,  g(b) \ne g(\overline{b})\,  z(b)
\end{align*}
So, for instance,  if we take  
\[
  g(b) =z(b)^{-1}, \ \  g(\overline{b})= g\,  z(b)^{-1} \mbox{ with }  g\ne e, 
\]
then $v(z)$ is not flat. Note that the above choice is always possible 
because $G$ is nontrivial by assumption. 
In conclusion we have shown the following 
\begin{teo}
\label{Cc:5}
Let $K$ be a pathwise connected but not  totally directed poset. 
Let $G$ be a nontrivial group. Then for any $1$--cocycle 
$z\in\Zrm^1(K,G)$ there are  connections in $\Urm^1(K,z)$ which are 
not flat. 
\end{teo}
Concerning central connections, in the case that $K$ and $G$ satisfy the
hypotheses of  the statement of  Theorem \ref{Cc:5}, and 
the centre of the group $G$ is nontrivial, then 
by using the above reasoning it is very easy to prove the existence
of nonflat central connections. 
\subsection{Holonomy and reduction of connections}
\label{Cd}
Consider a connection  $ u$ of  $\Urm^1(K,G)$. 
Fix  a base $0$--simplex $a_0$ and define 
\begin{equation}
\label{Cd:1}
 H_u(a_0)  \defi  \big\{  u(p)\in G \ |  \ p\in K(a_0)\big\} \ , 
\end{equation}
recalling that  $K(a_0)$ is the set of loops of $K$ 
with endpoint $a_0$. By the defining properties of connection
$1$--cochains it is very easy to see that  $H_u(a_0)$ is a
subgroup of $G$. Furthermore let 
\begin{equation}
\label{Cd:2}
 H_u^0(a_0)  \defi  \{  u(p)\in G \ |  \ p\in
 K(a_0), \ p\sim \sigma_0(a_0)\} \ ,  
\end{equation}
where $p\sim \sigma_0(a_0)$ means that $p$ is homotopic to the degenerate
$1$--simplex $\sigma_0(a_0)$. In this case, too, it is easy to see that 
$H_u^0(a_0)$ is a subgroup of $G$. Moreover, since 
$p*q*\overline{p}\sim\sigma_0(a_0)$ whenever 
$q,p\in K(a_0)$ and  $q\sim \sigma_0(a_0)$,  
$H^0_u(a_0)$ is a normal subgroup of
$H_u(a_0)$. $H_u(a_0)$
and $H^0_u(a_0)$ are called respectively 
\emph{the holonomy and the restricted holonomy group of $u$ 
         based on $a_0$}.\\ 
\indent As $K$ is pathwise connected, 
we have the following 
\begin{lemma}
\label{Cd:3}
Given  $u\in\Urm^1(K,G)$, let $\ga:G_1\to G$ 
be an injective homomorphism. The following assertions hold.\\  
(a) $H_u(a_0)$ and $H_u(a_1)$ are conjugate subgroups of $G$
    for any $a_0,a_1\in\Si_0(K)$. \\ 
(b) Given $u_1\in\Urm^1(K,G_1)$. If $\gamma\circ u_1$
    is equivalent to $u$, then 
           the holonomy groups 
     $H_{u_1}(a_0)$ and $H_u(a_0)$ are isomorphic.\\
The same assertions hold for the restricted holonomy groups. 
\end{lemma}
\begin{proof}
(a) Let $p$ be a path from $a_0$ to $a_1$. For any $g\in H_u(a_0)$, 
there is a loop $q\in K(a_0)$ such that $g=u(q)$.
Observe that $p*q*\overline{p}\in K(a_1)$, hence 
$u(p)\,  g\,  u(p)^{-1} =u(p*q*\overline{p})\in H_u(a_1)$.
By the symmetry of the reasoning, 
$H_u(a_0)\ni g\rightarrow  u(p)\,  g\,  u(p)^{-1}\in 
 H_u(a_1)$ is a group  isomorphism. (b) Let $u_1\in \Urm^1(K,G_1)$ 
and let $\fr\in(\ga\circ u_1,u)$. Since 
for any loop $p\in K(a_0)$, $\fr_{a_0}\, \ga\circ u_1(p) = 
 u(p)\,  \fr_{a_0}$, the map 
$H_{u_1}(a_0)\ni g\rightarrow 
\fr_{a_0}\,  \ga(g)\,  \fr^{-1}_{a_0}\in 
 H_u(a_0)$ is a group isomorphism.
\end{proof}
We now prove an analogue of the Ambrose-Singer theorem
for connections  of a poset. 
\begin{teo}
\label{Cd:4}
Let $u\in\Urm^1(K,G)$, $a_0\in\Si_0(K)$ and let $\io$ be 
the inclusion of $H_u(a_0)$ in $G$. 
Then there exists $u_1\in\Urm^1(K, H_u(a_0))$ such that 
$\io\circ u_1\cong u$.
\end{teo}
\begin{proof}
For any $0$--simplex $a$, let $p_a$ be a path from $a_0$ to
$a$.  Then define 
\begin{equation}
u_1(b) \defi u(\overline{p_{\partial_0b}}* b*p_{\partial_1b}),
\qquad b\in\Si_1(K).
\end{equation}
Note that $u_1(b)\in H_u(a_0)$ for any $1$--simplex $b$
because $\overline{p_{\partial_0b}}* b*p_{\partial_1b}\in K(a_0)$. 
Secondly, for any $1$--simplex $b$ we have 
\[
u_1(\overline{b}) = 
 u(\overline{p_{\partial_1b}}* \overline{b}*p_{\partial_0b}) =  
u(\overline{\overline{p_{\partial_0b}}* b* p_{\partial_1b}}) =
u_1(b)^{-1} \ .
\]
Thirdly, let $c\in \Si^{\inf}_2(K)$. Then 
\begin{align*}
u_1(\partial_0c)\,  u_1(\partial_2c) & =
  u(\overline{p_{\partial_{00}c}}*\partial_{0}c*p_{\partial_{10}c})\,  
   u(\overline{p_{\partial_{02}c}}*\partial_{2}c*p_{\partial_{12}c})
  \\
& =
  u(\overline{p_{\partial_{00}c}})\,  
  u(\partial_{0}c)\,  u(p_{\partial_{10}c})\,  
  u(\overline{p_{\partial_{02}c}})\,  
  u(\partial_{2}c)\,  u(p_{\partial_{12}c}) \\
& =
  u(\overline{p_{\partial_{01}c}})\,  
  u(\partial_{0}c)\,  
  u(\partial_{2}c)\,  u(p_{\partial_{11}c}) \\
& =
  u(\overline{p_{\partial_{01}c}})\,  
  u(\partial_{1}c)\,  u(p_{\partial_{11}c}) \\
& = u_1(\partial_{1}c) \ .
\end{align*}
Therefore we have that $u_1\in\Urm^1(K,H_u(a_0))$.
Finally, for any $0$--simplex $a$ let $\fr_a\defi u(p_a)$. 
Then for any $1$--simplex $b$ we have  
\begin{align*}
 \fr_{\partial_0b}\,  u_1(b) & =  
  u(p_{\partial_0b})\,  u(
  \overline{p_{\partial_0b}}*b*p_{\partial_1b}) \\
 & =  u(p_{\partial_0b})\,  u(\overline{p_{\partial_0b}})
      \,  u(b)\,  u(p_{\partial_1b}) 
 =  u(b)\,  u(p_{\partial_1b}) \\
 & =  u(b)\  \fr_{\partial_1b} \ ,
\end{align*}
namely $\fr\in(\io\circ u_1,u)$. 
Thus $\io\circ u_1\simeq u$.
\end{proof}
\section{Gauge transformations}
\label{D}
In the previous sections we have given several results to support 
the interpretation of $1$--cocycles of a poset as 
principal bundles over the poset. As the final issue of the present paper, 
we now introduce what we mean by the group of gauge transformations 
of a $1$--cocycle. \\[5pt] 
\indent Given a $1$--cocycle $z$ of $\Zrm^1(K,G)$,  define 
\begin{equation}
\label{D:1}
\mathcal{G}(z) \defi (z,z).
\end{equation}
An element of $\mathcal{G}(z)$ will be 
denoted by $\gr$. The composition law 
between morphisms of  $1$--cochains  endows 
$\mathcal{G}(z)$ with a structure of a group. The identity 
$\mathrm{e}$ of this group is given by $\mathrm{e}_a=e$ for any 
$0$--simplex $a$. The inverse $\g^{-1}$ of 
an element $\g\in\mathcal{G}(z)$ is obtained by composing $\g$ with the
inverse of $G$. We call $\mathcal{G}(z)$ the  
\emph{group of gauge transformations of $z$}.  
\begin{lemma}
\label{D:2}
If $z\in\Brm^1(K,G)$, then $\mathcal{G}(z)\cong G$.
\end{lemma} 
\begin{proof}
Observe that, since $K$ is connected, $\mathcal{G}(\imath)$ is the set of 
constant functions from $\Si_0(K)$ to $G$ and hence is isomorphic to $G$. 
As $z$ is a $1$--coboundary, 
it is equivalent to the trivial $1$--cocycle 
$\imath$, i.e.  there exists an $\fr\in(z,\imath)$. 
The mapping  
$\mathcal{G}(\imath)\ni\gr\mapsto \fr^{-1}\  \gr\
\fr\in\mathcal{G}(z)$ is a group isomorphism.
\end{proof}
As a consequence of this lemma and Proposition \ref{Bc:16},  
if the poset is simply connected then $\mathcal{G}(z)\cong G$
for any $1$--cocycle $z$. This is also the case when $G$ is Abelian.
\begin{lemma}
\label{D:3}
If $G$ is Abelian, then $\mathcal{G}(z)\cong  G$ for any 
$z\in\Zrm^1(K,G)$.
\end{lemma}
\begin{proof}
For any $\g\in\mathcal{G}(z)$ and for any $1$--simplex $b$ we have 
\[
\gr_{\partial_1b}\,z(b) = z(b)\, \gr_{\partial_0b} = 
\gr_{\partial_0b}\,z(b) \ . 
\]
Hence $\gr_{\partial_1b}=\gr_{\partial_0b}$ for any $1$--simplex $b$. 
Since $K$ is
pathwise connected, $\gr_{a} = g$ for any $0$--simplex $a$.
\end{proof}
Thus, for Abelian groups, the action of the group of 
gauge transformations is always \emph{global}, that is independent of 
the $0$--simplex.\smallskip

Given a $1$--cocycle $z\in\Zrm^1(K,G)$ consider the group 
$\mathcal{G}(z)$  of gauge transformations of $z$. 
For any $\ur\in\Urm^1(K,z)$ and 
$\g\in \mathcal{G}(z)$, define 
\begin{equation}
\label{D:4}
  \alpha_\g(u)(b)\defi \g_{\partial_0b}\,  u(b)\, \g_{\partial_1b}^{-1},
  \qquad b\in\Si_1(K).
\end{equation}
We have the following
\begin{prop}
\label{D:5}
Given $z\in\Zrm^1(K,G)$, the following assertions hold:\\
(a) given  $\g\in \mathcal{G}(z)$, then 
$\alpha_\g(u)\in \Urm^1(K,z)$ for any $u\in\Urm^1(K,z)$;\\
(b) The mapping 
\begin{equation}
\label{D:6}
\alpha:  \mathcal{G}(z)\times \Urm^1(K,z)\ni (\g,u) 
  \ \ \longrightarrow \ \ \alpha_\g(u) \in \Urm^1(K,z)
\end{equation}
defines a left action, not free,  of $\mathcal{G}(z)$ on
 $\Urm^1(K,z)$. 
\end{prop}
\begin{proof}
(a) Clearly $\alpha_\g(u)(\overline{b})= \alpha_\g(u)(b)^{-1}$ 
for any $1$--simplex $b$. Moreover, if $b\in\Si^{\inf}_1(K)$,  then  
$\alpha_\g(u)(b) = $ 
$\g_{\partial_0b}\,  u(b)\,  \g_{\partial_1b}^{-1}$
$= \g_{\partial_0b}\,  z(b)\,  \g_{\partial_1b}^{-1}$
$= z(b)$. This entails that $\alpha_\g(u)$ satisfies property 
(ii) of the definition of connections. Hence
$\alpha_\g(u)\in \Urm^1(K,z)$.
(b) Clearly, $\alpha$ is a left action that is not free, because
$z\in \Urm^1(K,z)$, hence $\alpha_\g(z) = z$ for any 
$\g\in\mathcal{G}(z)$. 
\end{proof}
\section{Conclusions and outlook}
We have developed a theory of bundles over posets from
a cohomological standpoint, the analogue of describing the usual
principal bundles in terms of their transition functions. In a
sequel, we will introduce principal bundles over posets and
their mappings directly and further develop such concepts
as connection, curvature, holonomy and transition function 
(we will also introduce concepts such as gauge group and gauge
transformation). Although all these concepts are familiar from
the usual theory of principal bundles, at this point it is
worth stressing some of the differences from that theory. As we
shall see in the sequel, the definition of principal bundle
involves bijections between different fibres satisfying a
$1$--cocycle identity. An important r$\hat o$le is played
by the simplicial set of inflationary simplices. All principal
bundles can be trivialized on the fundamental covering. Finally,
it should be stressed that the goal of these investigations is to
develop gauge theories in the framework of algebraic quantum
field theory. Our principal fibre bundles and the associated
vector bundles are envisaged stepping stones to the algebra of
observables. \\[8pt]

\noindent{\small \textbf{Acknowledgements.} 
 We wish to thank Ezio Vasselli for discussions on the topic 
 of this paper.}

\end{document}